\documentclass[12pt]{article}
\usepackage{latexsym}
\title{Riemannian manifolds not quasi-isometric to
leaves in codimension one foliations}
\author{Paul A. Schweitzer, S.J.\thanks{Supported in early stages of this
work by the CNPq, FAPERJ, PRONEX of the Ministry of Science and
Technology of Brazil, and the Clay Mathematics Institute. Address:
Depto. de Matem\'atica, Pontif\'\i cia Universidade Cat\'olica do
Rio de Janeiro (PUC-Rio), Rio de Janeiro, RJ 22453-900, Brasil;
email: paul37sj@gmail.com}}

\usepackage{amssymb}
\usepackage{amsmath}
\usepackage{graphics}
\usepackage[dvips]{color,graphicx}
\usepackage{psfrag}
\usepackage{float}

\newtheorem{thm1}{Theorem}[section]
\newtheorem{thm}{Theorem}
\newcommand{\bthm}{\begin{thm}}
\newcommand{\ethm}{\end{thm}}
\newtheorem{prop}[thm1]{Proposition}
\newcommand{\bpr}{\begin{prop}}
\newcommand{\epr}{\end{prop}}
\newtheorem{defn}[thm1]{Definition}
\newcommand{\bdf}{\begin{defn}}
\newcommand{\edf}{\end{defn}}
\newtheorem{lem}[thm1]{Lemma}
\newcommand{\blm}{\begin{lem}}
\newcommand{\elm}{\end{lem}}
\newtheorem{cor}[thm1]{Corollary}
\newcommand{\bcr}{\begin{cor}}
\newcommand{\ecr}{\end{cor}}
\newtheorem{conj}[thm1]{Conjecture}
\newcommand{\bcj}{\begin{conj}}
\newcommand{\ecj}{\end{conj}}
\newtheorem{rmk}[thm1]{Remark}
\newcommand{\brk}{\begin{rmk}}
\newcommand{\erk}{\end{rmk}}
\newtheorem{ex}[thm1]{Example}
\newcommand{\bex}{\begin{ex}}
\newcommand{\eex}{\end{ex}}

\newcommand{\F}{\mathcal{F}}
\newcommand{\T}{\mathcal{T}}
\newcommand{\B}{\mathcal{B}}
\newcommand{\R}{\mathcal{R}}
\newcommand{\V}{\mathcal{V}}
\newcommand{\Hor}{\mathcal{H}}

\newcommand{\sm}{\smallsetminus}

\newcommand{\Int}{{\rm Int}}
\newcommand{\qed}{\hfill$\Box$}

\begin{document}
\maketitle

\begin{abstract} Every open manifold $L$ of dimension
greater than one has complete Riemannian metrics $g$ with bounded
geometry such that $(L,g)$ is not quasi-isometric to a leaf of a
codimension one foliation of a closed manifold. Hence no
conditions on the local geometry of $(L,g)$ suffice to make it
quasi-isometric to a leaf of such a foliation. We introduce the
`bounded homology property', a semi-local property of $(L,g)$ that
is necessary for it to be a leaf in a compact manifold in
codimension one, up to quasi-isometry. An essential step involves
a partial generalization of the Novikov closed leaf theorem to
higher dimensions.

\medskip

\noindent{\bf keywords:} codimension one foliation, Reeb
component, non-leaf, geometry of leaves, bounded homology property

\medskip

\end{abstract}

CONTENTS

1. Introduction

2. Definitions and statements of results 

3. Leaves have the bounded homology property

4. Modifying the metric on a manifold by inserting balloons.

5. Proofs of two lemmas

6. Invariance under quasi-isometry

7. Novikov's Theorem for embedded $1$-connected vanishing cycles

References.

\section{Introduction}

The question of when an open (i.e. noncompact) connected manifold
can be realized up to diffeomorphism as a leaf in a foliation of a
compact differentiable manifold was first posed by Sondow
\cite{So} for surfaces in $3$-manifolds, and was solved positively
for all open surfaces by Cantwell and Conlon \cite{CC}. In the
opposite direction, Ghys \cite{G} and independently Inaba,
Nishimori, Takamura and Tsuchiya \cite{INTT} constructed an open
$3$-manifold (an infinite connected sum of lens spaces for all odd
primes) that cannot be a leaf in a foliation of a compact
$4$-manifold. Attie and Hurder \cite{AH} gave an uncountable
family of smooth simply connected $6$-dimensional manifolds, all
having the homotopy type of an infinite connected sum of copies of
$S^2\times S^4$, that are not diffeomorphic to leaves in a compact
$7$-manifold. It is still an open problem whether every smooth open
manifold of dimension greater than 2 is diffeomorphic to a leaf
of a codimension two (or higher) foliation.

In the related question of when an open Riemannian manifold can be
realized up to quasi-isometry as a leaf in a foliation of a
compact manifold, Attie and Hurder \cite{AH} also produced  an
uncountable family of quasi-isometry types of Riemannian metrics
on the $6$-manifold $S^3\times S^2\times {\mathbb R}$, each with
bounded geometry, which cannot be leaves in any codimension one
foliation of a compact $7$-manifold. On these and other
$6$-manifolds they also defined Riemannian metrics that have
positive `entropy', and hence cannot be leaves in any $C^{2,0}$
codimension one foliation, or in any $C^1$ foliation of arbitrary
codimension, on any compact manifold. Zeghib \cite{Z} adapted this
result to surfaces with exponential growth. Attie and Hurder prove
their various results using the bounded Pontryagin classes defined
by Januskiewicz \cite{J}, who had already used them to construct
open manifolds that could not be leaves in certain compact
manifolds, extending earlier results of Phillips and Sullivan
\cite{PS}.

The Attie-Hurder results extend to codimension one foliations of
dimensions greater than 6, but their Question 2 asks for examples
in the lower dimensions $3,4$, and $5$. In this paper we respond
to this question by showing that {\em every} open manifold of
dimension at least $3$ admits complete Riemannian metrics of
bounded geometry, and of every possible growth type, that are not
quasi-isometric to leaves of codimension one foliations of compact
manifolds. The same result for surfaces had already been proven in
\cite{Sc1} (where the present paper was announced). Thus no set of
local bounds on the geometry of an open Riemannian $p$-manifold
$L$ with $p\geq 2$ can be sufficient to guarantee that it be
quasi-isometric to a leaf of a $C^{2,0}$ codimension one foliation
of a closed $(p+1)$-manifold.

We define a $C^{2,0}$ foliation to be one in which the leaves are
smooth of class $C^2$, and their $C^2$ differentiable structure
varies continuously in the transverse direction along the leaves
of a transverse foliation $\T$ (Definition \ref{C20}). Attie and
Hurder \cite{AH} call such foliations $C^0$. A bound on the
geometry of $L$ is {\bf local} if it only depends on the
Riemannian geometry of the $\epsilon$-balls around the points of
$L$, for some constant $\epsilon>0$. For example, bounds on the
various curvatures and on the injectivity radius of $L$ are local
geometric bounds.

We define two Riemannian manifolds to be {\bf quasi-isometric} if
there is a diffeomorphism $f$ from one to the other such that both
$f$ and $f^{-1}$ produce only bounded distortion of the metrics,
up to a constant $D$, as formulated precisely in Definition
\ref{qi} below and in \cite{AH}. It is easy to see that a leaf in
a foliation of a smooth compact manifold $M$, with any Riemannian
metric induced by a metric on $M$, must have {\bf bounded
geometry}, in the sense that the sectional curvature is uniformly
bounded above and below, and the injectivity radius has a uniform
positive lower bound. Because of the term $D$ in Definition \ref{qi},
quasi-isometry does not always preserve this property, but in this
paper we consider only Riemannian metrics with bounded geometry.
Then the following theorem, which is our main result, involves a
restriction on the global geometry of an open manifold for it to
be quasi-isometric to a leaf.

\bthm Every connected non-compact smooth $p$-manifold $L$ with
$p\geq 2$ possesses $C^\infty$ complete Riemannian metrics $g$
with bounded geometry that are not quasi-isometric to any leaf of
a codimension one $C^{2,0}$ foliation on any compact
differentiable $(p+1)$-manifold.

Furthermore $g$ can be chosen such that no end is quasi-isometric
to an end of a leaf of such a foliation, and also to have any
growth type compatible with bounded geometry. Hence there are
uncountably many quasi-isometry classes of such metrics $g$ on
every such manifold $L$. \label{mainthm}\ethm

\begin{figure}[H]
\centering
      \includegraphics*[width=.5\linewidth]{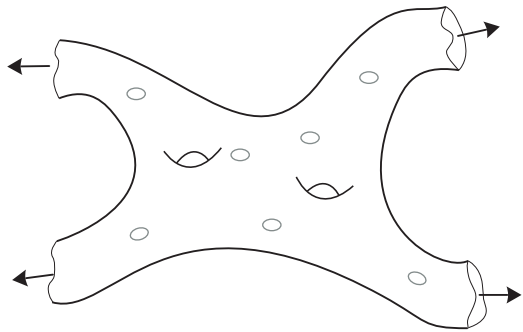}
\caption{The manifold $L$ with the original
metric.}\label{originalsurf}
\end{figure}

When $p=2$, Theorem \ref{mainthm} was proven in \cite{Sc1}, so in
this paper we shall assume that $p\geq 3$. The proof for $p\geq 3$
depends on Theorem \ref{thm2} below, which states that the leaves
of a $C^{2,0}$ codimension one foliation of a compact manifold $M$
of dimension greater than $3$, with Riemannian metrics that vary
continuously in the transverse direction, have a certain {\bf
bounded homology property}, a property which we define for open
Riemannian manifolds of dimension at least $3$. In Theorem
\ref{thmexs} we show how to modify a given complete Riemannian
metric on a connected smooth open manifold, without changing its
growth type, so that the new metric does not have the bounded
homology property. The construction involves introducing spherical
`balloons' of arbitrarily large size, but with `necks' of
uniformly bounded size, as in  \cite{Sc1}. (See Figures
\ref{originalsurf} and 2.)


\begin{figure}[H]\label{withballoons}
\centering
      \includegraphics*[width=.5\linewidth]{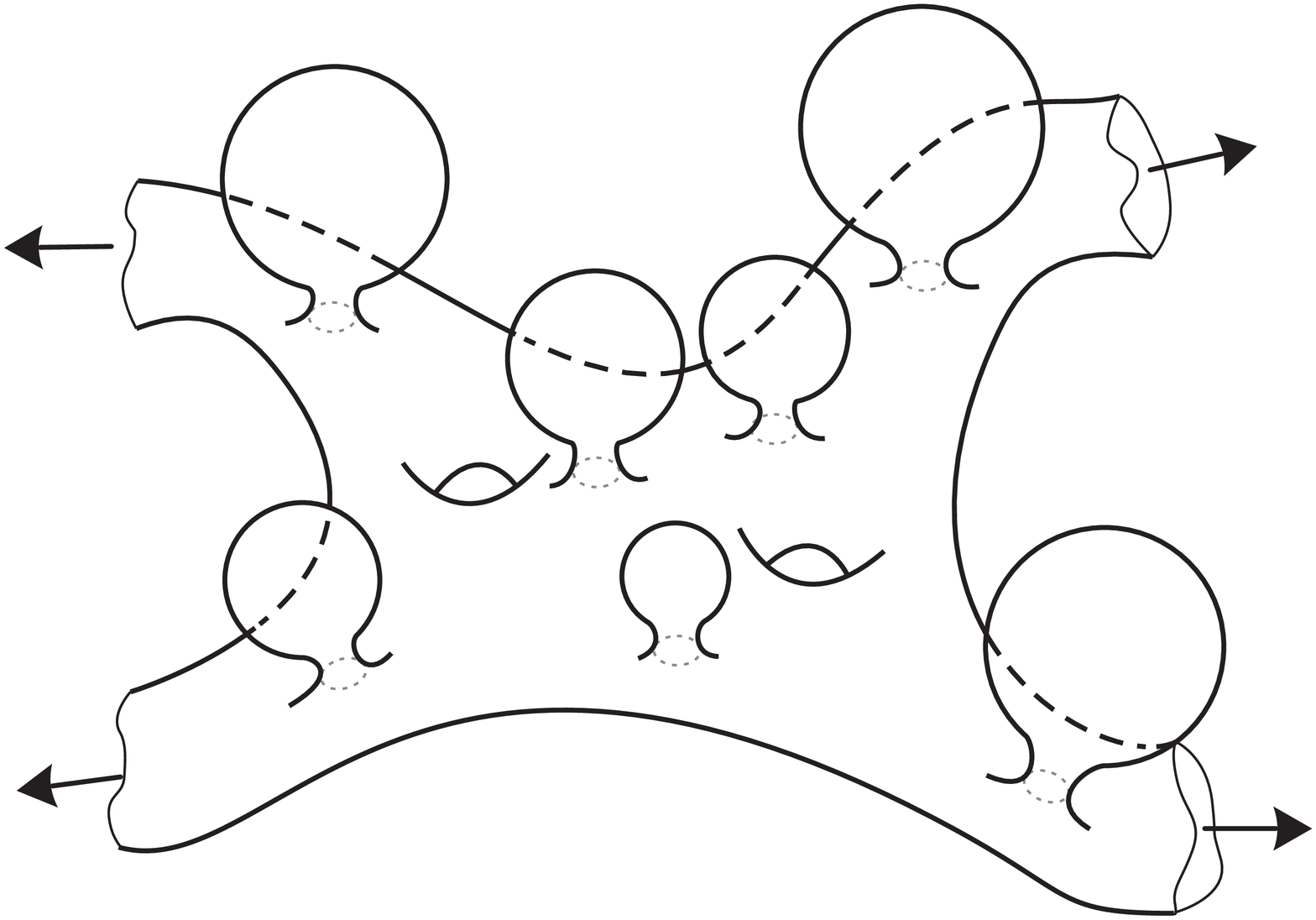}
\caption{The manifold $L$ with ``balloons''.}
\end{figure}


The proof of Theorem \ref{thm2} (in Section \ref{pf}) involves a
Finiteness Lemma (Lemma \ref{fntlemma}, proved in Section
\ref{fntsection}) and a partial extension to higher dimensions of
Novikov's theorem that vanishing cycles only occur on the boundary
of Reeb components (Theorem \ref{Nov-thm}, proved in Section
\ref{Novikov-gen}). Definitions and statements of results are
given in Section \ref{defs} and the construction of Riemannian
manifolds that cannot be leaves is given in Section \ref{exs}.

The first three theorems of this paper constitute an extension to
higher dimensions of three similar theorems for open surfaces in
\cite{Sc1}. That paper, using a different bounded {\em homotopy}
property for surfaces involving contractible loops rather than
bounding submanifolds, showed that all open surfaces have complete
Riemannian metrics with bounded geometry that cannot be leaves in
$C^{2,0}$ foliations of compact $3$-manifolds.

It is a pleasure to thank my student Fabio Silva de Souza for
preparing the figures.

\section{Definitions and statements of results}\label{defs} 

In this section we give several definitions leading up to the
definition of the bounded homology property, and then we state
Theorem \ref{thm2} (which states that, under certain hypotheses,
leaves have this property), Theorem \ref{thmexs}, and Theorem
\ref{Nov-thm} (a partial generalization of Novikov's Theorem on
the existence of Reeb components). In this paper, all manifolds
are smooth of class $C^2$, except that in the proof of Theorem
\ref{Nov-thm} it suffices to assume that $M$ and $\F$ are $C^0$.

\bdf A diffeomorphism $f:L\rightarrow L'$ between two Riemannian
manifolds $L$ and $L'$ is a {\bf quasi-isometry} if there exist
constants $C, D>0$ such that the distance functions $d$ and $d'$
on $L$ and $L'$ satisfy
$$C^{-1}d'(f(x),f(y))-D\leq d(x,y) \leq Cd'(f(x),f(y))+D$$
for all points $x,y\in L$. When such a diffeomorphism exists we
say that $L$ and $L'$ are {\bf quasi-isometric}. \label{qi}\edf
For example, any diffeomorphism between compact smooth Riemannian
manifolds is a quasi-isometry. Note that quasi-isometry is an
equivalence relation.

Let $L$ be a $p$-dimensional Riemannian manifold, $S$ a subset of
$L$, and $\beta$ a positive number. The open $\beta$-ball
$V_\beta(x)$ at a point $x\in L$ is defined to be the set of
points on $L$ whose distance from $x$ is less than $\beta$. Note that
$V_\beta(x)$ may fail to be a topological ball, if $\beta$ is
greater than the injectivity radius at $x$.

\bdf The $\beta$-{\bf volume} ${\rm Vol}_\beta(S)={\rm
Vol}_{\beta,L}(S)\in {\mathbb N} \cup\{\infty\}$ of $S$ on $L$ is
the smallest integer $K$ such that $S$ can be covered by $K$ open
$\beta$-balls in $L$, or $\infty$ if no such finite number exists.
\edf

The $\beta$-{\bf volume} of $S$ depends on its embedding in $L$,
but the ambient manifold $L$ will be clear from the context, so it
will not be made explicit in the notation.

\bdf Let $C$ be a compact Riemannian $p$-manifold with boundary
$B$ and let $\beta$ be a constant greater than $0$. Then we define
the {\bf Morse $\beta$-volume} $M(C,\beta)$ to be the smallest
positive integer $M$ for which there is a Morse function
$f:C\rightarrow [0,\infty)$ such that $f(B)=0$ and each level set
$f^{-1}(t)$ for $t\in [0,\infty)$ has ${\rm
Vol}_\beta(f^{-1}(t))\leq M$ on $C$. \label{Morsedeltadef} (See
Figure 3, where various level sets are shown.) \edf


\begin{figure}[H]\label{Morsevol}
\centering

\includegraphics*[width=.7\linewidth]{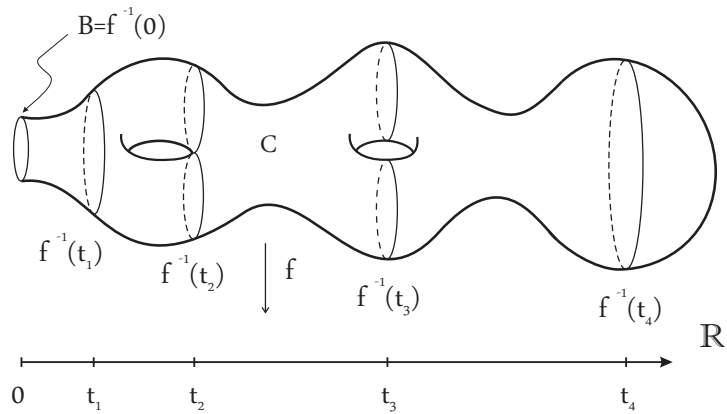}
\caption{Morse $\beta$-volume of a compact manifold $C$ with
boundary $B$.}
\end{figure}
This definition makes sense since it is evident that the
$\beta$-volumes of the level sets of a fixed Morse function $f$ on
a compact manifold are uniformly bounded, so that $M(C,\beta)$ is
finite. For a set $S$ contained in a Riemannian manifold $L$, we
define the open {\bf $\beta$-neighborhood} of $S$, $V_\beta(S)$,
to be the set of all points in $L$ whose distance from $S$,
measured along geodesics in $L$, is less than $\beta$. We can now
define the bounded homology property, the fundamental tool used in
this paper.

\bdf A Riemannian $p$-manifold $L$ ($p\geq 3$) has the {\bf
bounded homology property} if there exists a constant $\beta_0\geq
0$ such that for every pair of constants $\beta>\beta_0$ and $k>0$
there is an integer $K=K(\beta, k)$ such that if
\begin{enumerate}
\item $B\subset L$ is a $1$-connected smooth closed $(p-1)$-submanifold
embedded in $L$,
\item there is a closed neighborhood $V$ of $B$ that fibers over
$B$ as a smooth tubular neighborhood and contains the
$\beta$-neighborhood $V_\beta(B)$ of $B$ in $L$,
\item $B$ has  $\beta$-volume ${\rm Vol}_\beta(B)\leq k$ on $L$, and
\item $B$ bounds a compact 1-connected region $C$
on $L$, \end{enumerate} then the Morse $\beta$-volume of $C$
satisfies $M(C,\beta)\leq K$. \label{BHP} \edf {\bf Comments.} A
neighborhood $V$ of $B$ is a {\bf tubular neighborhood} of $B$ if
there is a smooth retraction $V\rightarrow B$ which is a (locally
trivial) fibration. In particular, it is well-known that for a
smooth compact submanifold $B$ of a Riemannian manifold
$V_\beta(B)$ is a tubular neighborhood of $B$ for every
sufficiently small positive number $\beta$. On the other hand,
because of the term $D$ in Definition \ref{qi}, we must require
$\beta$ to be `sufficiently large', i.e., greater than some given
$\beta_0$, in order for the bounded homology property to be
invariant under quasi-isometry (see Proposition \ref{2BHP}). We
also require that $B$ and $C$ be $1$-connected (i.e., connected
and simply connected) so that when they are subsets of a leaf $L$
of a foliation $\F$, they can be lifted in the transverse
direction to leaves of $\F$ near to $L$. Recall that for any $C^0$
codimension one foliation $\F$ there always exists a
$1$-dimensional foliation $\T$ topologically transverse to $\F$ \cite{HH2},
\cite{Si3}.

\bdf A codimension one foliation $\F$ of a smooth manifold $M$ is a
$C^{2,0}$-foliation if the leaves are smooth submanifolds of class $C^2$ and
the transverse foliation $\T$ can be chosen so that the
$C^2$-structures on the leaves are preserved by the local
homeomorphisms obtained by lifting open sets from one leaf to
another along leaves of $\T$. \label{C20}\edf

\bthm Let $\F$ be a $C^{2,0}$ codimension one foliation of
dimension $p\geq 3$ on a closed smooth manifold $M$ and let $\T$
be a foliation transverse to $\F$ chosen as in Definition
\ref{C20}. Then the leaves of $\F$ have the bounded homology
property with $\beta_0=0$, and for every $k>0$ and $\beta>0$, the
same integer $K=K(k,\beta)>0$ can be chosen for all the leaves of
$\F$.\label{thm2} \ethm

\bcr If $\F$ is a $C^{2,0}$ codimension one foliation of dimension
$p\geq 3$ of a closed smooth manifold $M$, then the leaves have
the bounded homology property for any Riemannian metric induced on
the leaves by a Riemannian metric on $M$. \ecr

The proof of Theorem \ref{thm2} is given in Section \ref{pf} using
results proven in Sections \ref{fntsection} and \ref{Novikov-gen}. In
Section \ref{exs} we show that every smooth open manifold of
dimension at least $3$ admits complete Riemannian metrics with
bounded geometry such that the bounded homology property does not
hold, thus giving the following result.

\bthm Let $L$ be a complete open connected Riemannian $p$-manifold
of dimension $p\geq 3$ whose metric $g_0$ has bounded geometry.
Then $L$ has other complete Riemannian metrics $g$ with bounded
geometry and the same growth type as $g_0$ that do not possess the
bounded homology property, and such that no end of $L$ has the
bounded homology property. Furthermore there are such metrics with
uncountably many distinct growth types, and hence in uncountably
many distinct quasi-isometry classes. \label{thmexs} \ethm
We recall that
the {\bf growth function} $f :[0,\infty)\rightarrow[0,\infty)$ of
a connected Riemannian $p$-manifold $(L,g)$ with basepoint $x_0\in
L$ is defined to be $f(r)={\rm Vol}(B(x_0,r))$, the
$p$-dimensional volume of the ball of radius $r$ centered at
$x_0$. Given two increasing continuous functions $f_1, f_2
:[0,\infty)\rightarrow[0,\infty)$, we say that $f_1$ has {\bf
growth type} less than or equal to that of $f_2$ (denoted
$f_1\preceq f_2$) if there exist constants $A,B,C>0$ such that for
all $r\in [0,\infty)$, $f_1(r)\leq Af_2(Br+C)$ \cite{H}. They have
the same growth type if $f_1\preceq f_2$ and $f_2\preceq f_1$. We
write $f_1\prec f_2$ if $f_1\preceq f_2$ but it is false that
$f_2\preceq f_1$. For example, $I\prec {\rm exp}$ where $I$ is the
identity function $I(r)=r$ and ${\rm exp}(r)=e^r$ on $[0,\infty)$.
The growth type of $(L,g)$ (i.e., of its growth function) is
clearly invariant under quasi-isometry and change of the
basepoint. Then the following observation proved in Section
\ref{exs} establishes the last assertions of Theorems
\ref{mainthm} and \ref{thmexs}, since the growth types of the
functions $f_k(r)=r^k$ for every $k>1$ are distinct and can all be
realized.

\bpr Let $f:[0,\infty)\rightarrow [0,\infty)$ be any increasing
continuous function of growth type greater than linear and at most
exponential, i.e. $I\prec f\preceq {\rm exp}$. Then every smooth
connected open manifold $L$ of dimension at least two admits a
complete Riemannian metric with bounded geometry whose growth type
is the growth type of $f$. \label{growth}\epr

Theorem \ref{mainthm} follows immediately from Theorems
\ref{thmexs} and \ref{thm2} with the following fact.

\bpr If $L$ and $L'$ are quasi-isometric complete Riemannian
manifolds with bounded geometry and $L$ has the bounded homology
property, then so does $L'$. \label{2BHP}\epr

To prove this proposition we need the following result. Both will
be proven in Section \ref{qipfsect}.

\bpr Let $L$ be a complete Riemannian manifold with bounded
sectional curvature. Then given constants $0<a<b$, there exists an
integer $n>0$ such that every open $b$-ball $V_b(x)$ on $L$
centered at $x\in L$ can be covered by at most $n$ open $a$-balls.
\label{2balls}\epr

As mentioned above, $V_b(x)$ may fail to be a topological ball.

\medskip

The proof of Theorem \ref{thm2} will use a weak generalization of
the second half of Novikov's theorem on the existence of Reeb
components and a Corollary, which we shall state after the
following definition.

\bdf A compact $(p+1)$-dimensional manifold with a codimension one
foliation $\R$ is a (generalized) {\bf Reeb component} if the
boundary $\partial R$ is a nonempty finite union of leaves, the
interior ${\rm Int}(R)$ fibers over the circle with the leaves as
fibers, and there is a transverse orientation pointing inwards
along all the components of $\partial R$ . \label{Reebcomp}\edf

It is clear that the boundary leaves are compact and the interior
leaves non-compact. In this paper, we shall usually consider Reeb
components with connected boundary $\partial R$, and then the
existence of the transverse orientation pointing inwards along
$\partial R$ is automatic.

Now let $\F$ be a $p$-dimensional foliation of a compact
$(p+1)$-dimensional manifold $M$ and let $B$ be a compact
connected $(p-1)$-dimensional manifold. The horizontal foliation
$\Hor$ of $B\times [0,1]$ is given by the leaves $B\times \{t\}$
for $t\in [0,1]$. Consider a foliated map
$$h: (B\times [0,1],\Hor)\rightarrow (M,\F)$$  and suppose that
$h_0: B\rightarrow L_0$ is an embedding,
where for all $t$ $L_t$ is the leaf containing $h(B\times\{t\})$ and
$h_t: B\rightarrow L_t$ is the map defined
$h_t(b)=h(b,t)$. Note that now we are not supposing any differentiability,
but there does exist a $1$-dimensional foliation $\T$ of $M$ topologically
transverse to $\F$.
Then we have the following weak generalization of the second half
of Novikov's theorem.

\bthm {\bf An Extension of Novikov's Theorem.}
Suppose that one of the following conditions holds for every
$t>0$ sufficiently close to $0$ but does not hold for $t=0$:
\begin{enumerate}
\item $B_t=h_t(B)$ is
the boundary of a compact $1$-connected region $C_t\subset L_t$;
\item $B_t=h_t(B)$ is
the boundary of a compact region $C_t\subset L_t$;
\item $B$ is oriented and $0=h_{t*}([B])\in H_{p-1}(L_t)$ (where $[B]$ is the
fundamental homology class of $B$); or
\item $0=h_{t*}([B])\in H_{p-1}(L_t;{\mathbb Z}_2)$ (where $[B]$ is the
fundamental homology class of $B$ with coefficients modulo $2$).
\end{enumerate}
Then the leaf $L_0$ is the boundary of a Reeb component $R$ whose
interior ${\rm Int}(R)$ is the union of the leaves $L_t$ for which
$t>0$. (See Figure 4.) \label{Nov-thm}\ethm

\begin{figure}[H]\label{classicalRC}
\centering
\includegraphics*[width=.5\linewidth]{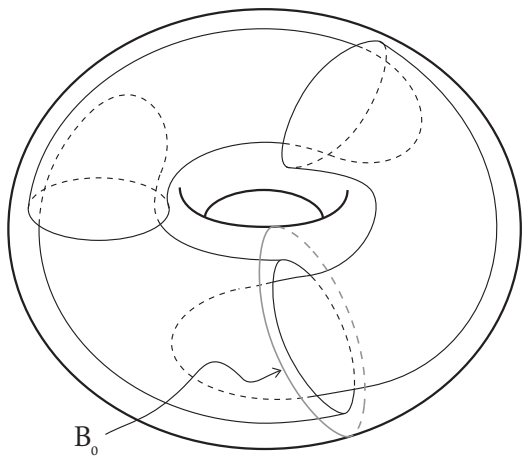}
\caption{The classical Reeb component with vanishing cycle $B_0$.}
\end{figure}

\bcr Assume the same hypotheses and also that $\F$ is a $C^{2,0}$
foliation. Then for any $\beta>0$, and for any choice of
Riemannian metrics on the leaves in the Reeb component that vary
continuously along the transverse foliation $\T$, there is a
constant $K>0$ such that the Morse $\beta$-volume
$M(C_t,\beta)\leq K$ for all $t\in (0,1]$. \label{MorseReeb} (See
Figure 5.) \ecr
\begin{figure}[H]\label{MorsevolinRC}
\centering
\includegraphics*[width=.6\linewidth]{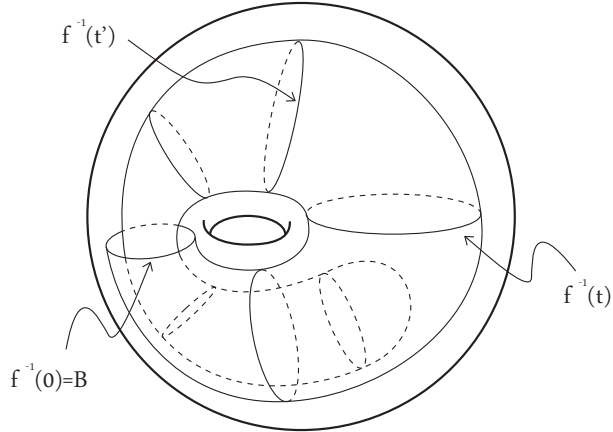}
\caption{The Morse $\beta$-volume of a set $C$ in a leaf of a Reeb
component.}
\end{figure}

For the results of this paper, we only use condition (1) of
Theorem \ref{Nov-thm}, but the other three conditions are
mentioned since they are of some interest and involve little extra work.
In the Corollary we assume $\F$ to be smooth of class $C^{2,0}$ in
order for the leaves to admit Riemannian metrics and Morse functions on the regions $C_t$.
Note that as $t$ approaches $0$, the Riemannian volume of $C_t$
and also its $\beta$-volume ${\rm Vol}_\beta(C_t)$ both tend to
infinity, as is easily seen from the structure of the Reeb
component. That is the reason for using the Morse $\beta$-volume
$M(C_t,\beta)$, which is uniformly bounded according to the
Corollary, instead of the other two volumes, which are not, in the
Definition \ref{BHP} of the bounded homology property.

The proofs of Theorem \ref{Nov-thm} and its Corollary are given in
Section \ref{Novikov-gen}. We observe that the proof in this case,
in which $B$ is an {\em embedded} ``vanishing cycle'', is much
easier than in the general case when $B$ is only assumed to be
{\em immersed}. (See Haefliger \cite{H2} or Camacho-Lins Neto
\cite{CL} for good expositions of the proof of Novikov's original
theorem that treat the problem of double points of the immersion
clearly.) In the course of the proof of Theorem \ref{Nov-thm}, we
find a construction of Reeb components that is used in the proof
of the Corollary. In \cite{AHS}, where the structure of
generalized Reeb components is studied in detail, it is shown that
this construction actually produces all Reeb components, up to
foliated homeomorphism.

\section{Leaves have the bounded homology property}\label{pf} 

In this section we prove Theorem \ref{thm2} using results that
will be proven in Sections \ref{fntsection} (the Finiteness Lemma
\ref{fntlemma}) and \ref{Novikov-gen} (the extension of Novikov's
Theorem \ref{Nov-thm} and Corollary \ref{MorseReeb}). The idea of
the proof is to consider a $(p-1)$-dimensional submanifold $B$
that satisfies the hypotheses of Definition \ref{BHP} and
approximate it by a subcomplex $X$ of a triangulation on the same
leaf. Then we show that there is a finite set of ``transverse
families'' of subcomplexes that suffice in this process (the
Finiteness Lemma \ref{fntlemma}). Finally we show that there is a
constant $K$ that is an upper bound for all the resulting Morse
$\beta$-volumes.

Let $\F$ be a codimension one foliation of a compact manifold $M$
of dimension at least $4$, so that the leaf dimension $p\geq 3$,
and suppose that constants $k>0$ and $\beta>0$ are given. Suppose
that the leaves are smooth (of class $C^2$) and their $C^2$
differentiable structure varies continuously in the transverse
direction along leaves of a fixed transverse foliation $\T$, so
that $\F$ is $C^{2,0}$, as in Definition \ref{C20}. As a
consequence, it is possible to choose Riemannian metrics on the
leaves of $\F$ which vary continuously in the transverse direction
along $\T$, and we fix such Riemannian metrics. Thus the
hypotheses of Theorem \ref{thm2} are satisfied. All these
structures remain fixed throughout this section.

As usual, a {\bf region} in a smooth manifold $L$ is a compact
connected submanifold with smooth boundary that is the closure of
an open set in $L$. If $L$ has a Riemannian metric, then the
diameter of a smooth submanifold $S$ in $L$ (possibly with
boundary and corners, such as a simplex) is the supremum of
distances between pairs of points of $S$, as measured along
geodesics in $S$. The {\bf mesh} of a triangulation of a region is
the maximum of the diameters of its simplexes. A triangulation is
an {\bf $\epsilon$-triangulation} if its mesh is less than a
positive number $\epsilon$. All of the triangulations we shall
consider will be $\beta'$-triangulations with $\beta'=\beta/4$.

\vskip 0.5cm \noindent {\bf The Simplicial Approximation Process
(SAP).} Let $\B$ be the set of all submanifolds $B$ of leaves of
$\F$ that satisfy the conditions (1), (2), and (3) of Definition
\ref{BHP}. Given $B\in \B$ on a leaf $L$ of $\F$, choose a smooth
$\beta'$-triangulation of a region in $L$ that contains the
tubular neighborhood $V$ given by condition (2), and give $B$ a
smooth triangulation that is sufficiently fine so that the open
star of each vertex $v$ of $B$ lies in the open star of some
vertex, say $g(v)$, of the triangulation of $L$. Then the function
$g$ taking the vertices of $B$ to those of $L$ extends to a
simplicial map $\bar g: B\rightarrow L$, which is a simplicial
approximation to the inclusion $B\hookrightarrow L$. Note that the
image $X_B=\bar g(B)$ is a $(p-1)$-dimensional complex contained
in $V_{\beta'}(B)$. (See Figure 6, which shows $B$ and its
simplicial approximation $g(B)$ in the tubular neighborhood $V$
which is contained in the region $\Omega$.)

We say that $X_B$ is obtained from $B$ by the {\bf SAP}, the
simplicial approximation process.
\begin{figure}[H]\label{simplapprox}
\centering
\includegraphics*[width=.6\linewidth]{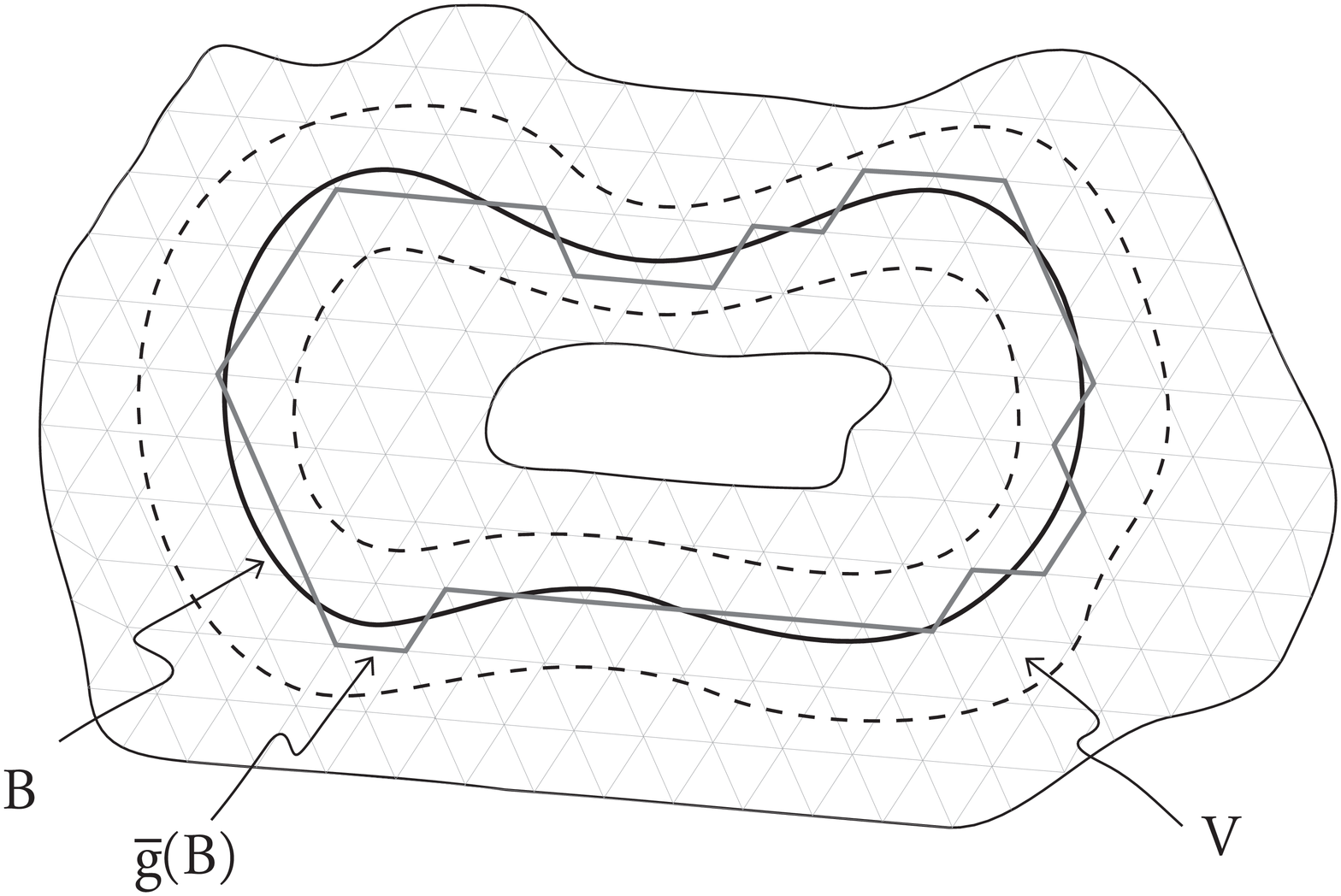}
\caption{$B$ and its simplicial approximation $\bar g(B)$ in the
neighborhood $V$.}
\end{figure}
\vskip 0.5cm

\noindent {\bf Transverse families and the Finiteness Lemma.}
Recall that a product manifold $X\times Y$ has two product
foliations $\Hor$ and $\V$, the horizontal and vertical
foliations, given respectively by the leaves $X\times \{y\}$ for
$y\in Y$ and $\{x\}\times Y$ for $x\in X$. A map $f: (X\times
Y,\Hor, \V)\rightarrow (M,\F,\T)$ is {\bf bifoliated} if it takes
leaves of $\Hor$ and $\V$ into leaves of $\F$ and $\T$,
respectively.

\bdf A triple $(X,\Omega, f)$ will be called a {\bf transverse
family} if $\Omega$ is a smoothly triangulated compact $p$-manifold
with boundary, $X\subset \Omega$ is a $(p-1)$-dimensional
subcomplex, and
$$f: (\Omega\times [0,1],\Hor,\V)\rightarrow (M,\F,\T)$$ is a
smooth bifoliated embedding, such that if we set
$\Omega_t=f(\Omega\times \{t\})$ and $X_t=f(X\times \{t\})$, then,
for every $t\in [0,1]$, relative to the metric on the leaf $L_t$
that contains $\Omega_t$ and $X_t$,
\begin{enumerate}
\item $\Omega_t$
contains $V_{2\beta'}(X_t)$, the $2\beta'$-neighborhood of $X_t$
(with $\beta'=\beta/4$), and
\item the triangulation induced on
$\Omega_t$ by the triangulation on $\Omega$ is a smooth
$\beta'$-triangulation.\end{enumerate} \label{trcomplex}\edf
 We say that a transverse family
$(X,\Omega,f)$ is {\bf good} if there is a submanifold  $B\subset
L_0$ with $B\in\B$ and a smooth bifoliated embedding
$$h: (B\times [0,1],\Hor,\V)\rightarrow (M,\F,\T)$$
such that for every $t\in [0,1]$ the submanifold $B_t=h(B\times
\{t\})$ is contained in $L_t$, $B_0=B$, and $X_t$ can be obtained
from $B_t$ by the SAP described above.

The following Finiteness Lemma will be proven in Section
\ref{fntsection}.

\blm {\bf (Finiteness Lemma).} There is a finite set of good
transverse families  $(X_i, \Omega_i, f_i), 1\leq i\leq \ell$,
such that, for each $B\in \B$, the complex $X_B$ obtained from $B$
by the SAP can be chosen to be $X_{i,t}=f_i(X_i\times \{t\})$ for
some $i\in \{1,\dots,\ell\}$ and $t\in [0,1]$.
\label{fntlemma}\elm

For $1\leq i\leq \ell$, let $\B_i$ be the set of $B\in\B$ for
which there exists some $t\in [0,1]$ for such that the SAP can
yield the complex $X_B=X_{i,t}$.

\bpr For each $i, 1\leq i\leq \ell,$ there exists an integer $K_i$
such that, for each $B\in \B_i$ that satisfies condition (4) of
Definition \ref{BHP}, the region $C$ with $\partial C=B$ given by
condition (4) has Morse $\beta$-volume $M(C,\beta)$ less than or
equal to $K_i$. \label{boundx}\epr Supposing Lemma \ref{fntlemma}
and Proposition \ref{boundx}, we can now prove Theorem \ref{thm2}.

\medskip

\noindent {\bf Proof of Theorem \ref{thm2}.}  By the Finiteness
Lemma, we obtain a finite set of good transverse families $(X_i,
\Omega_i, f_i), 1\leq i\leq \ell$. The Proposition gives a common
upper bound $K_i$ for the Morse $\beta$-volumes of the regions $C$
corresponding to all the submanifolds $B\in \B_i$. Since every
$B\in \B$ belongs to some $\B_i$, $K=K(k,\beta)={\rm max}_{1\leq
i\leq \ell}K_i$ is an upper bound for the Morse $\beta$-volumes
$M(C,\beta)$ of the regions $C$ corresponding to all submanifolds
$B$ that satisfy the hypotheses of the Theorem, as claimed. \qed

\medskip

\noindent {\bf Proof of Proposition \ref{boundx}.} First of all,
we observe that if some $B\in \B$ bounds $1$-connected regions on
both sides in its leaf $L$, then the leaf $L$ is compact, and by
the Van Kampen Theorem it is also $1$-connected. Then by the Reeb
Stability Theorem for codimension one, all the leaves are compact
and the foliation $\F$ fibers over the circle with the leaves as
fibers. Thus there is a common upper bound for the $\beta$-volumes
of the leaves, and hence for all the regions $C$, so that the
conclusion of the Proposition holds. (In fact, the conclusion of
Theorem \ref{thm2} also holds.) Hence we may assume that each
$B\in \B$ bounds a compact $1$-connected region on at most one
side.

Now fix an index $i$ as in the Finiteness Lemma \ref{fntlemma} and
consider the good transverse family $(X_i, \Omega_i,f_i)$, which
for simplicity we denote by $(X,\Omega,f)$ (without indicating the
index). Set \[J=\{t\in [0,1]\ |\ \exists{\rm\ a\
1\!\!-\!connected\ region}\ C_t\subset L_t {\rm\ such\ that}\
\partial C_t=B_t\}.\]
 Thus the indices $t\in J$ are those for
which $B_t=B_{i,t}$ satisfies condition (4)  of Definition
\ref{BHP}. It is clear that $J$ is open in the interval $[0,1]$,
since the region $C_t$ is $1$-connected by hypothesis, and
therefore lifts along $\T$ to nearby leaves.

Observe that no connected component of $J$ can be an open interval
$(a,b)$, for then the submanifolds $B_a$ and $B_b$ do not bound
$1$-connected regions (since $a,b\notin J$), but the leaves $L_t$
for $t$ in the interval $(a,b)$ do. Then, by the extension of
Novikov's Theorem, Theorem \ref{Nov-thm}, the leaves $L_a$ and
$L_b$ bound Reeb components, with interior leaves $L_t$ for $t\in
(a,1]$ in the first Reeb component, and $L_t$ for $t\in [0,b)$ in
the second one. But then $L_a$ is compact as the boundary of the
first Reeb component, but it is non-compact, since it is an
interior leaf of the second one, a contradiction.

Consequently the only possibilities for the connected components
of $J$ are $[0,b)$, $(a,1]$, and $[0,1]$. In the first two cases,
$L_b$ or $L_a$ is the boundary of a Reeb component with the leaves
in the indicated interval as interior leaves, and then Corollary
\ref{MorseReeb} gives a common upper bound for the Morse
$\beta$-volumes of the corresponding regions $C_t$. Note that if
for some $t_1\in J,\ M(C_{t_1},\beta)\leq K$ for a certain integer
$K$, then $M(C_t,\beta)\leq K$ for all $t$ in a small neighborhood
of $t_1$, since the compact regions $C_t$ vary continuously and
$C_{t_1}$ is covered by $K$ open $\beta$-balls. Thus the Morse
$\beta$-volumes $M(C_t,\beta)$ are locally bounded, and if
$J=[0,1]$, the compactness of $[0,1]$ gives the desired common
upper bound. If $J=\emptyset$, then there is nothing to be proven.
Putting all these cases together, we get a common upper bound
$K'_i$ for all the $C_t$ with boundary $B_t\in \B_i$, $t\in
J\subset [0,1]$.

To complete the proof of the proposition we need the following
lemma, whose proof is given in Section \ref{fntsection}.

\blm There is an integer $K_0$ that depends on $k, \beta,$ and the
upper bound on the scalar curvature of the leaves of $\F$, such
that if the same complex $X_t$ can be obtained by the SAP from
each of two submanifolds $B, B' \in \B$, and $B$ satisfies
condition (4) of Definition \ref{BHP}, i.e., there is a compact
$1$-connected region $C$ on the leaf $L$ containing $B$ such that
$\partial C=B$, then $B'$ also satisfies condition (4), so there
is a compact $1$-connected region $C'\subset L$ whose boundary is
$B'$, and we have
\begin{equation} |M(C',\beta) - M(C,\beta)| \leq K_0.
\end{equation} \label{B1B2} \elm
Every $B\in \B_i$ yields a complex $X_{i,t}$ under the SAP for
some $t\in [0,1]$, and the submanifold $B_t=B_{i,t}$ yields the
same complex, so the lemma shows that $K_i=K'_i + K_0$ is an upper
bound for the Morse $\beta$-volumes of all the regions $C$ with
$\partial C=B\in \B_i$, as claimed. \qed

\section{Modifying the metric on a manifold by inserting balloons.} \label{exs}

In this section we prove Theorem \ref{thmexs} by showing how the
Riemannian metric on any complete open Riemannian manifold $L$ of
dimension at least three with bounded geometry can be modified
without changing the growth type, so that $L$ with the new metric
does not have the bounded homology property, and hence cannot be a
leaf in a $C^{2,0}$ codimension one foliation of a closed
manifold. The construction is an obvious adaptation of the
construction for surfaces in Section 2 of \cite{Sc1}, which we
follow closely. We insert $p$-dimensional ``balloons'' of
unbounded size with ``necks'' of uniformly bounded size into $L$.
The balloons are widely spaced so that the original growth type of
$L$ does not change.

\medskip
\noindent{\bf Proof of Theorem \ref{thmexs}.} Let $L$ be a smooth
connected noncompact $p$-dimen\-sion\-al manifold ($p\geq 3$) and
let $g_0$ be a complete Riemannian metric on $L$ with globally
bounded sectional curvature, with injectivity radius greater than
some small constant $d>0$, and with a given growth type. Hence the
exponential map ${\rm exp}_x: B_x(0,d)\rightarrow L$ is injective
for every point $x$ in $L$, where $B_x(0,r)$ denotes the open ball
of radius $r$ centered at the origin in the tangent space at $x$.
Consequently ${\rm exp}_x$ is a diffeomorphism from $B_x(0,d)$
onto the open ball $V_d(x)$ in $L$.

We fix a basepoint $x_0\in L$ and consider a sequence of positive
numbers $d_n, n = 1,\dots,$ with $d_n+2d<d_{n+1}$ for each $n$.
Choose a sequence of points $x_1,x_2,\dots$ such that
$d(x_0,x_n)=d_n$ and consequently the balls $V_d(x_n)$ are
disjoint. Since the metric $g_0$ has bounded geometry, there is no
loss of generality in supposing, as we do, that the balls
$V_d(x_n)$ are isometric to the balls $V_d(0)$ in Euclidean
$p$-space ${\mathbb R}^p$. In fact, it suffices to modify $g_0$
slightly (perhaps replacing $d$ by a smaller radius) with care to
preserve the growth type and globally bounded geometry.

Choose an increasing sequence of positive numbers $r_n$ such that
$r_n\rightarrow \infty$ as $n\rightarrow \infty$. Let $S^p(r_n)$
be the the sphere of radius $r_n$ in ${\mathbb R}^{p+1}$ centered
at the origin and let $S=(0,\dots,0,-r_n)$ be its basepoint at the
south pole. Choose a diffeomorphism $\phi: V_d(x_n)\rightarrow
S^p(r_n)\sm \{S\}$ by setting
$$\phi({\rm exp}_{x_n}(tv)) = {\rm exp}_S((d-t)h(v))$$
for $d/2\leq t <d$ and every unit tangent vector $v$ to $L$ at
$x_n$, where $h$ is a linear isometry from the tangent space to
$L$ at $x_n$ to the tangent space to $S^p(r_n)$ at $S$, and then
extending $\phi$ as a diffeomorphism over the complementary disk
$V_{d/2}(x_n)$. Define a new metric $g$ on $V_d(x_n)$ by
interpolating between the given metric $g_0$ and the metric $g_1$
obtained as the pullback under $\phi$ of the round metric on
$S^p(r_n)$, so that $g$ coincides with $g_0$ near the boundary of
$V_d(x_n)$ and with $g_1$ on $V_{d/2}(x_n)$. This defines the new
``balloon'' metric on the balls $V_d(x_n)$, and outside these
balls $g$ is defined to coincide with the original metric $g_0$,
as shown in Figure 2. The uniform manner of carrying out this
construction as $n$ varies ensures that $(L,g)$ has globally
bounded geometry. Note that the metric $g$ on the closed ball
$\overline{V_{d/2}(x_n)}$ of radius $d/2$ in the original metric
$g_0$ is now the round metric of $S^p(r_n)$ from which a small
neighborhood of $S$ has been removed.

We claim that the new metric $g$ does not have the bounded
homology property. In fact, suppose that $\beta_0\geq 0$ is given
and fix a number $\beta>\beta_0$. Take $n_0$ sufficiently large so
that $r_{n_0}>2\beta+2d$, and consider only the balloons for
$n\geq n_0$. Let $B_n$ be the $(p-1)$-sphere on this balloon with
radius $\beta+d$ centered at the south pole $S$ in the original
metric on the sphere. Note that $B_n$ has $(p-1)$-volume less than
$a_{p-1}(\beta+d)^{p-1}$, where $a_m$ denotes the $m$-volume of
the unit $m$-sphere $S^m$ in ${\mathbb R}^{m+1}$, and $B_n$ is
$1$-connected since $p\geq 3$. Choose $k$ to be sufficiently large
so that $S^{p-1}$ can be covered by $k$ open balls of radius
$\beta$; it follows that the same is true for $B_n$. Furthermore,
it is clear that the closed $\beta$-neighborhood $V =
\overline{V_\beta(B_n)}$ of $B_n$ in the new metric is a tubular
neighborhood fibered over $B_n$ that contains $V_\beta(B_n)$. The
closed complementary component $C_n$ of $B_n$ that contains the
north pole $N=(0,\dots,r_n)$ is a compact $1$-connected region on
$L$ with boundary $\partial C_n=B_n$. Hence the four conditions of
Definition \ref{BHP} are satisfied by $B_n$, but we shall see that
the Morse volumes of the $C_n$ are unbounded.
\begin{figure}[H]\label{Morsevolballoon}
\centering
\includegraphics*[width=.5\linewidth]{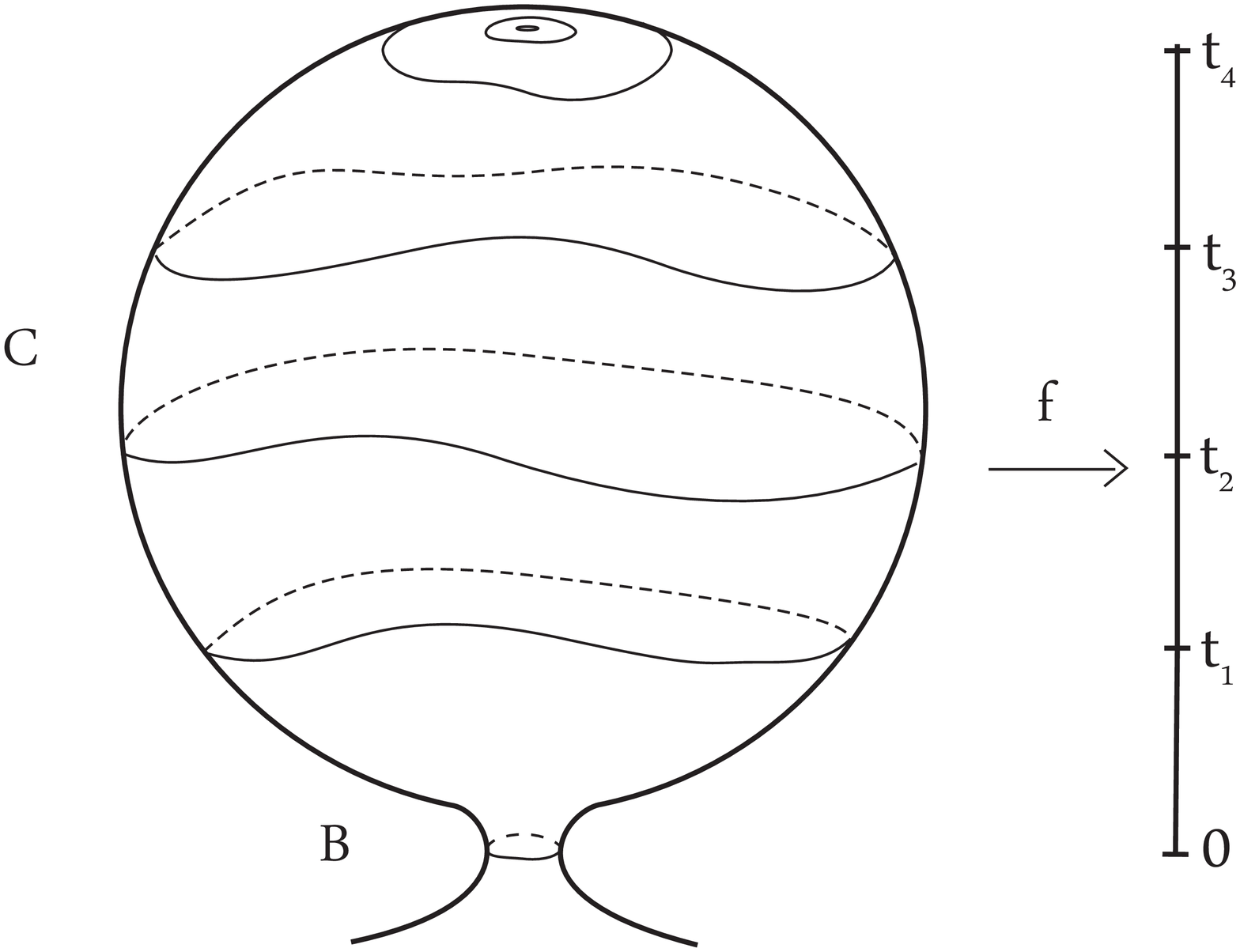}
\caption{The Morse $\beta$-volume of a balloon $C$.}
\end{figure}

Let $f_n: C_n\rightarrow [0,\infty)$ be a Morse function with
$f_n(B_n)=0$. As $t$ increases from $0$, there is a $1$-parameter
family of closed complementary regions
$A_n(t)=f_n^{-1}([t,\infty))$ of $f_n^{-1}(t)$, beginning with
$A_n(0)=C_n$ and ending with $A_n(t)=\emptyset$ for sufficiently
large values of $t$. Note that $A_n(0)$ covers more than half of
the balloon and the $p$-volume of the sets $A_n(t)$ varies
continuously, so for some value of $t$, say $t_n'$, $A_n(t_n')$
will have $p$-volume $a_p{r_n}^p/2$, i.e. $A_n(t_n')$ will have
half the $p$-volume of the round sphere of radius $r_n$, which is
$a_p{r_n} ^p$. Then, by the isoperimetric inequality on the
sphere, the boundary $f_n^{-1}(t_n')=\partial A_n(t_n')$ must have
$(p-1)$-volume at least equal to $a_{p-1}{r_n}^{p-1}$, the
$(p-1)$-volume of the equator. (See Figure 7.) As $n$ increases
these $(p-1)$-volumes tend to infinity, so each $\beta$-volume
${\rm Vol}_\beta(f_n^{-1}(t_n'))$ is greater than some constant
$M_n$ such that ${\rm lim}_{n\rightarrow\infty}M_n=\infty$. It
follows that the Morse $\beta$-volume of $C_n$ satisfies
$M(C_n,\beta)>M_n$, so $(L,g)$ does not have the bounded homology
property.

The points $d_n$ can be chosen so that no end of $L$ has the
bounded homology property. In fact, if the points $x_i$ have been
chosen for $i\leq n$, let $k$ be the number of connected
components of the complement of $V_{d_n}(x_0)$ with noncompact
closure, and choose $x_{n+1},\dots,x_{n+k}$ so that one of them is
in each of these components. Continue inductively, repeating this
procedure with $n+k$ in place of $n$. This guarantees that each
end will contain an infinite number of the points $x_n$ and
therefore does not have the bounded homology property.

We can choose the sequence $\{d_n\}$ to grow sufficiently fast,
relative to the sequence $\{r_n\}$, so that the growth functions
$f_0$ of $g_0$ and $f$ of $g$ satisfy $f_0(r)\leq f(r)\leq
2f_0(r)$. Hence $g$ will have the same growth type as $g_0$, so
the process of inserting balloons can be carried out so as to
preserve the growth type.

The last conclusion, that there are uncountably many
quasi-isometry classes of Riemannian metrics on $L$ that do not
have the bounded homology property, follows from Proposition
\ref{growth}, which realizes the distinct growth types $x\mapsto
x^k$ for every $k>1$. \qed

\medskip
\noindent {\bf Proof of Proposition \ref{growth}.} We outline the
argument, leaving precise details to be filled in by the reader.
Given a smooth connected non-compact manifold $L$ of dimension
$p\geq 2$ with basepoint $x_0$ and any continuous increasing
function $f:[0,\infty)\rightarrow [0,\infty)$ such that $I\prec
f\preceq {\rm exp}$, we shall find a metric $g_2$ on $L$ such that
the growth function is equivalent to $f$.

Let $q:L\rightarrow [0,\infty)$ be a proper smooth Morse function
such that $q^{-1}(0)=\{x_0\}$. Choose Riemannian metrics on each
level set $q^{-1}(t)$ such that they vary smoothly in $t$ away
from the singularities of $q$ and the geometry of these level sets
is uniformly bounded. We construct a Riemannian metric $g_1$ on
$L$ so that, away from the singularities of $q$, the level sets
are totally geodesic, the gradient flow of $\phi\circ f$ (for a
diffeomorphism $\phi:[0,\infty)\rightarrow [0,\infty)$ to be
chosen later) is orthogonal to the level sets, and the distance
between the level sets $q^{-1}(s)$ and $q^{-1}(t)$ is
$|\phi(t)-\phi(s)|$; the metric must be adjusted in small
neighborhoods of the singularities. If the function $\phi$ grows
sufficiently rapidly, the level sets will be spread far apart in
comparison with the growth of their volumes, and the growth
function $f_1$ of the metric $g_1$ will satisfy $f_1\preceq f$.
Here we must use the hypothesis that $f$ grows more than linearly,
since the requirement of globally bounded geometry and the
topology of the level sets $q^{-1}(t)$ may make it impossible to
have a uniform upper bound on their volumes.

Now if $f_1\prec f$ we can increase the growth type by inserting
sufficiently many large balloons in balls $V_d(x_n)$ (and possibly
in other disjoint balls $V_d(y)$ that are closer together,
including balls on balloons already inserted) to change $g_1$ to a
new metric $g_2$ so that the new growth function $f_2$ will have
the same growth type as the given function $f$. Note that it is
possible to get the exponential growth type by inserting so many
balloons that their number grows exponentially as a function of
the distance from the basepoint $x_0$, but since the geometry is
required to have bounded sectional curvature, it is impossible to
get any greater growth type. Clearly any intermediate growth type
can be realized by choosing an appropriate distribution of
inserted balloons. \qed

\brk It is an open question whether or not leaves of foliations of
codimension greater than one on closed manifolds have the bounded
homology property. Tsuboi's construction of a codimension two
foliation whose $2$-dimensional leaves do not have the bounded
{\em homotopy} property (see the last section of \cite{Sc1}) uses
the fact that certain loops on the leaves and are of unbounded
length since they are connected and have unbounded diameter. That
construction cannot be adapted to give leaves that do not have the
bounded {\em homology} property, since the level sets of Morse
functions need not be connected. In fact, the obvious adaptation
of the $2$-dimensional construction to higher dimensions gives
leaves that do possess the bounded homology property. \erk

\section{Proof of two lemmas}\label{fntsection} 

In this Section we prove Lemmas \ref{fntlemma} and \ref{B1B2}. We
assume all the conditions indicated in the first two paragraphs of
Section \ref{pf}.

\medskip
\noindent {\bf Proof of the Finiteness Lemma \ref{fntlemma}.}
Given a point $x\in M$ in a leaf $L_x$, let $D_x =
\overline{B(x,\epsilon)}$ be the closed $p$-dimensional ball
centered at $x$ of radius $\epsilon$ for some positive $\epsilon$
less than the injectivity radius at $x$, so that $D_x$ is actually
homeomorphic to a closed ball. We may lift $D_x$ along $\T$ to
disks $D_y$ on nearby leaves $L_y$ for $y\in J$, where $J$ is an
open interval embedded in a leaf of $\T$ and containing $x$, thus
obtaining a bifoliated map
$$h:(D_x\times J,\Hor,\V)\rightarrow (M,\F,\T)$$ such that
$h(x,y)=y\in L_y$. Fix a larger smooth compact region $E\subset
L_x$ that contains $V_d(D_x)$, where $d=(2k+1)\beta$, and give $E$
a smooth $\beta'$-triangulation, where as before $\beta'=\beta/4$.

If $B\in \B$ meets $D_x$, then $V_\beta(B)\subset E$ since $B$ is
connected and has $\beta$-volume less than or equal to $k$ and
therefore diameter less than $2k\beta$. Then the SAP applied to
$B$ will yield a $(p-1)$-dimensional subcomplex $X_B\subset E$.
Let $\Omega_B$ be a smooth compact region containing
$\overline{V_{2\beta'}(B)}$ in its interior and contained in
$V_{3\beta'}(B)$, and give $\Omega_B$ a smooth
$\beta'$-triangulation that agrees with the triangulation on $E$
for all the simplices contained in $V_{2\beta'}(B)$. The closed
tubular neighborhood $V$ of $B$ given by condition (2) of
Definition \ref{BHP} satisfies $\Omega_B\subset
V_{3\beta'}(B)\subset V$ and is $1$-connected, since by hypothesis
$B$ is. Therefore $\Omega_B$ can be lifted along $\T$ to
$\Omega_{B,y}$ on every leaf $L_y$ sufficiently close to $L_x$.
Hence we get a bifoliated embedding

$$h_B:(\Omega_B\times J',\Hor,\V)\rightarrow (M,\F,\T)$$ such that
$h_B(x',y)=y\in L_y$ for every point $y$ in a sufficiently small
open subinterval $J'\subset J$ containing $x$ and for some $x'\in
B\cap D_x$.

Now since $E$ is compact there are only finitely many
$(p-1)$-dimensional subcomplexes of $E$; let $X_1,\dots X_m$ be
those that are obtained by the SAP as $X_B$ for some $B\in \B$
that meets $D_x$. Choose a submanifold $B_i$ for each $X_i$, so
that $X_i=X_{B_i}$, ie., $X_i$ is obtained from $B_i$ by the SAP.
We may choose $J'$ small enough so that for every $B\in
\{B_1,\dots,B_m\}$ and $y\in J'$, the lifted triangulation on
$\Omega_{B,y}=h_B(\Omega_B \times \{y\})$ has mesh less than
$\beta'$ and $\Omega_{B,y}$ contains the the
$2\beta'$-neighborhood $V_{2\beta'}(X_{B,y})$ of the lifted
complex  $X_{B,y}=h_B(X_B \times \{y\})$ of $X_B$. Thus the
conditions (1) and (2) of Definition \ref{trcomplex} are satisfied
for $\Omega_{B,y}$ and $X_{B,y}$ on their leaves, for every $y\in
J'$. Now if we let $g: [0,1]\rightarrow J'$ be an embedding such
that $x= g(t)$ for some $t\in (0,1)$, the triple $(X_B, \Omega_B,
f_B)$, where $f_B(t,y)=h_B(g(t),y)$, will be a good transverse
family. Thus we get $m$ good transverse families corresponding to
$D_x$. Every $B\in \B$ that meets $D_y$ for $y\in g([0,1])$ will
yield an $X_{B,y}$ under the SAP.

For each point $x\in M$ we get an open set $h({\rm Int}(D_x)
\times g(0,1))$ containing $x$. Such open sets cover $M$, so
finitely many of them suffice to cover $M$, say for $x_1, \dots,
x_n$. The union of the sets of triples $(X_B, \Omega_B, f_B)$
where $B$ varies over all the $B_i's$ for all the points
$x_1,\dots, x_n$, will be a finite set of good transverse
families. It satisfies the conclusion of the Finiteness Lemma
since every $B\in \B$ will meet one of the sets  $h({\rm Int}(D_x)
\times g(0,1))$, for some $x\in\{x_1,\dots,x_n\}$, and so will
yield one of the complexes $X_{B_i}$ for that $x$. \qed
\medskip

\noindent{\bf Proof of Lemma \ref{B1B2}.} We suppose the
hypotheses of the Lemma. Observe that the SAP moves points of each
of $B$ and $B'$ a distance less than $\beta'$, since the
triangulation on the leaf has mesh less than $\beta'$. Hence
$B\subset V_{\beta'}(X_t)$ and $X_t\subset V_{\beta'}(B)$, and
similarly for $B'$. Thus $B'\subset V_{2\beta'}(B)$ and $B\subset
V_{2\beta'}(B')$. Now each of the two connected submanifolds $B$
and $B'$ separates the leaf $L$ into two connected components. For
$B$ they are the interior of $C$ and $L\sm C$. Since $B'\subset
V_{2\beta'}(B)$, one of the two connected components of $L\sm B'$
must be contained in $V_{2\beta'}(C)=C \cup V_{2\beta'}(B)$; call
its closure $C'$. Since $B\subset V_{2\beta'}(B')$, we must have
$C\subset V_{2\beta'}(C')= C'\cup V_{2\beta'}(B')$.

Clearly $\partial C'=B'$ and $C'$ is connected. We must show that
its fundamental group is trivial. By hypothesis, there is a closed
tubular neighborhood $V'$ of $B'$ containing $V_\beta(B')$. Set
$C'_-=C'\sm {\rm Int}(V')$ and $C'_+=C'\cup V'$, so that
$C'_-\subset C'\subset C´_+$. These inclusions induce isomorphisms
on the fundamental groups, since $C'$ is the union of $C'_-$ with
a collar neighborhood of its boundary that is homeomorphic to
$\partial C'\times I$ while $C'_+$ is the union of $C'$ with a
collar neighborhood of its boundary also homeomorphic to $\partial
C'\times I$. Since $B\subset V_{2\beta'}(B')$, $(C'\sm
V_{2\beta'}(B'))\cap B=\emptyset$, so $C'\sm
V_{2\beta'}(B')\subset C$, and we have $C'_-\subset C'\sm
V_{2\beta'}(B')\subset C\subset C'\cup V_{2\beta'}(B')\subset
C'_+$. Thus the inclusion $C'_-\subset C'_+$, which induces an
isomorphism of fundamental groups, factors through $C$, which is
$1$-connected by hypothesis, so the homeomorphic sets $C'_-, C'$,
and $C'_+$ are also $1$-connected.

Finally, we must show that there exists a constant $K_0$ such that
$M(C',\beta) \leq M(C,\beta) +K_0$, which by symmetry will
establish the second conclusion. Since the foliated manifold $M$
is compact, there is an upper bound for the sectional curvature of
the leaves of $\F$. By Proposition \ref{2balls} there is a
constant $n_0$, depending only on $\F$, $k$, and $\beta$, such
that every open $2\beta$-ball $V_{2\beta}(x)$ in a leaf of $\F$
can be covered by $n_0$ open $\beta$-balls in the leaf, and we set
$K_0=kn_0$. Now let $f: C\rightarrow [0,\infty)$ be a Morse
function on $C$ such that $f(B)=0$ and every level set $f^{-1}(t)$
has $\beta$-volume ${\rm Vol}_{\beta}(f^{-1}(t))\leq M(C,\beta)$.
Since $C\sm V_{2\beta'}(B)\subset {\rm Int}(C')$, we may extend
the Morse function $f+1$ restricted to $C\sm V_{2\beta'}(B)$ to a
Morse function $f': C'\rightarrow [0,\infty)$ such that
$f'^{-1}(0)=B'$. By hypothesis there are $k$ points
$y_1,\dots,y_k\in L$ such that the union of the $\beta$-balls
$V_\beta(y_j)$ covers $B'$, and then the union of the
$2\beta$-balls $V_{2\beta}(y_j)$ covers $V_{2\beta'}(B')$. Each of
these $2\beta$-balls can be covered by at most $n_0$
$\beta$-balls, so the $\beta$-volume of $V_{2\beta'}(B')$ is at
most $K_0=kn_0$. Finally, each level set $f'^{-1}(t)$ is contained
in the union $f^{-1}(t+1) \cup V_{2\beta'}(B')$, whose
$\beta$-volume is at most ${\rm Vol}_\beta(f^{-1}(t+1))  +{\rm
Vol}_\beta(V_{2\beta'}(B'))\leq M(C,\beta) + K_0$, showing that
$M(C',\beta) \leq M(C,\beta) +K_0$, as claimed. \qed

\medskip

\section{Invariance under quasi-isometry}\label{qipfsect} 

The goal of this section is to prove Proposition \ref{2BHP}, which
states that the bounded homology property is invariant under
quasi-isometry. We shall use Proposition \ref{2balls} and the
following result, which is certainly well known.

\bpr Let $L$ be a complete Riemannian manifold with sectional
curvature between $-c$ and $c$ for a constant $c\geq 1$ and let
$v_1,v_2\in S\subset T_xL$ be points on the unit sphere in the
tangent plane at a point $x\in L$. Then the distance between the
points ${\rm exp}_xtv_1$ and ${\rm exp}_xtv_2$ on $L$ is at most
$e^{ct}$ times the distance between $v_1$ and $v_2$ on $S$. \epr

\noindent{\bf Proof.} If the distance on the unit sphere $S
\subset T_xL$ from $v_1$ to $v_2$ is $d$, then there exists a
geodesic $v: [0,d]\rightarrow S$ from $v_1$ to $v_2$ parametrized
by arclength, so that $|v'(s)|\equiv 1$. Define
$f:[0,\infty)\times [0,d]\rightarrow L$ to be $f(t,s)={\rm
exp}_xtv(s)$, where ${\rm exp}_x$ is the exponential map. Then for
each $s$, $J_s(t)=\frac{\partial f}{\partial s}(t,s)$ is a Jacobi
vector field (see do Carmo \cite{dC}) and satisfies the Jacobi
equation
$$ \frac{D^2J}{dt^2}+R(\gamma'_s(t),J_s(t))\gamma'_s(t)=0$$
where $R$ is the curvature operator, $\frac{D}{dt}$ is the
covariant derivative, and $\gamma_s(t)={\rm exp}_xtv(s)$ is the
geodesic starting at $x$ with initial velocity
$\gamma_s'(0)=v(s)$. Fix a parallel orthonormal frame
$e_1,\dots,e_p$ along each $\gamma_s$, expand $J_s(t)$ as
$$J_s(t)=\sum_{i=1}^p f_{s,i}(t)e_i(t),$$ and set $a_{s,ij}=\ <R(\gamma_s'(t),
e_i(t)) \gamma_s'(t), e_j(t)>,$ where $<\cdot,\cdot>$ denotes the
Riemannian metric. Then the Jacobi equation translates to the
system of $p$ ordinary differential equations
$f_{s,j}''(t)=-\sum_i a_{s,ij}(t)f_{s,i}(t)$ (see \cite{dC}) or a
system of $2p$ first order differential equations
$Y'(t)=B(t)Y(t)$, where $Y(t)=(y_1(t),\dots,y_{2p}(t))$ with
$y_i(t)=f_{s,i}(t)$ and $y_{i+p}(t)=f'_{s,i}(t)$ for $1\leq i\leq
p$, and $B(t)$ is a $2p\times 2p$ matrix of the form
$$B(t)=\left[\begin{array}{cc}0 & I\\
-A(t) & 0\\ \end{array}\right]$$ where $A(t)=[a_{s,ij}(t)]$ and
$I$ is the identity matrix. Now $A(t)$ is symmetric  (by the
symmetry properties of the curvature $R$), and so we can
diagonalize it to $\bar A(t)$ in a new orthonormal frame $\bar
e_1,\dots,\bar e_n$. Since $R$ is multi-linear, each diagonal
coefficient of the new diagonal matrix $\bar A(t)$ will be
$$\bar a_{s,ii}=\ <R(\gamma_s'(t), \bar e_i(t)) \gamma_s'(t), \bar
e_i(t)>$$ which is just the sectional curvature in the plane of
$\gamma_s'(t)$ and $\bar e_i(t)$, so that $|\bar a_{s,ii}|\leq c$.
Hence the matrix $A(t)$ has norm $|A(t)|=|\bar A(t)|\leq c$ and so
$B(t)$ also has norm $|B(t)|\leq c$ (since $|I|= 1\leq c$). Note
that $J_s(0)=0$ and $J_s'(0)=\frac{\partial^2 f}{\partial
t\partial s}(0,s)=v'(s)$ with $|v'(s)|=1$, so $|Y(0)|=1$. Then
Theorem 1.5.1 of \cite{Zt} shows that $|J_s(t)|\leq |Y(t)|\leq
\exp(\int_0^t|B(r)|dr)\leq e^{ct}$. We have shown that the
velocity of the curve $s\mapsto {\rm exp}_xtv(s)$ is at most
$e^{ct}$ times that of the curve $v(s)$ on $S$, so the distance
between the points ${\rm exp}_xtv_1$ and ${\rm exp}_xtv_2$ on $L$
is at most $e^{ct}d$. \qed

\medskip

\noindent{\bf Proof of Proposition \ref{2balls}.} Let $L$ be a
complete Riemannian $p$-dimensional manifold with sectional
curvature between $-c$ and $c$ for some constant $c\geq 1$ and
suppose that constants $0<a<b$ are given. Take points $0\leq
t_1<t_2<\dots<t_r\leq b$ so that every $t\in [0,b]$ lies within a
distance less than $a/2$ of one of the points $t_j$. Let
$c_0=2e^{cb}$ and choose a set $\{v_1,\dots, v_m\}\subset S$ so
that every point of the unit sphere $S\subset T_xL$ is at most at
a distance $a/c_0$ on $S$ from one of the points $v_i$. Then every
point $y={\rm exp}_xt'v$ in the open ball $B(x,b)$ of radius $b$
centered at $x$ (for some $v\in S$ and $t'\in [0,b]$) lies within
a distance less than $e^{ct'}a/c_0\leq a/2$ from a point ${\rm
exp}_x(t'v_i)$ on a geodesic $t \mapsto{\rm exp}_x(tv_i)$, and
this point is at a distance at most $a/2$ from one of the points
${\rm exp}_x(t_jv_i)$. Therefore $y$ lies in one of the open balls
$B({\rm exp}_x(t_jv_i),a)$. Thus the open balls of radius $a$
centered at the $mr$ points ${\rm exp}_x(t_jv_i), i=1,\dots, m,
j=1,\dots, r$ cover the ball $B(x,b)$. The same argument applies
with the same values $m$ and $r$ around any other point $x'\in L$,
so every open $b$-ball on $L$ can be covered by at most $n=mr$
$a$-balls, as claimed. \qed

\medskip

\noindent {\bf Proof of Proposition \ref{2BHP}.} Let $h:
L'\rightarrow L$ be a quasi-isometry of two Riemannian manifolds
with bounded geometry, as in Definition \ref{qi}, so that there
are constants $C_0\geq 1$ and $D>0$ such that the quasi-isometry
inequalities
$$C_0^{-1}d(h(x'),h(y'))-D\leq d'(x',y') \leq C_0d(h(x'),h(y'))+D$$ hold for
all $x',y'\in L'$. Suppose that $L$ has the bounded homology
property, so there is a fixed $\beta_0\geq 0$ such that for all
$\beta>\beta_0$ and $k>0$, there is an integer $K=K(k,\beta)$ such
that any $B=\partial C\subset L$ satisfying the four conditions of
Definition \ref{BHP} must have $M(C,\beta)<K$. Let
$\beta'_0=C_0\beta_0 + D$ and suppose that numbers $k'>0$ and
$\beta'>\beta'_0$ are given. We must find an integer $K' =
K'(k',\beta')$ satisfying the bounded homology property on $L'$.

Set $\beta=C_0^{-1}(\beta'-D)>C_0^{-1}(\beta_0'-D)=\beta_0$,
$k=n_1k'$ and $K'=K$, where $n_1$ is a constant given by
Proposition \ref{2balls} such that on $L$ every ball of radius
$C_0(\beta'+D)$ is covered by at most $n_1$ balls of radius
$\beta$.

Note that if $V_r(x)$ and $V'_r(x')$ are the open $r$-balls on $L$
and $L'$ centered at $x\in L$ and $x'\in L'$, then by the
quasi-isometry inequalities
$$h(V'_{\beta'}(x'))\subset V_{C_0(\beta'+D)}(h(x'))$$ and
$$h^{-1}(V_{\beta}(h(x')))\subset V'_{\beta'}(x')$$
since $d'(x',y')<\beta'$ implies that
$d(h(x'),h(y'))<C_0(\beta'+D)$ and $d(h(x'),h(y'))<\beta$ implies
that $d'(x',y')<C_0\beta+D=\beta'$.

Consider any $(p-1)$-submanifold $B'$ of $L'$ satisfying the four
conditions of Definition \ref{BHP} for the constants $k'$ and
$\beta'$, so that there exists a tubular neighborhood $V'$ of $B'$
containing the $\beta'$-neighborhood $V'_{\beta'}(B')$ of $B'$,
$B'$ has $\beta'$-volume ${\rm Vol}_{\beta'}(B')\leq k'$ on $L'$,
and there exists a compact $1$-connected region $C'$ in $L'$ with
$\partial C'=B'$. We shall show that $M(C',\beta')\leq K'$.

Set $B=h(B'), V=h(V'),$ and $C=h(C')$. Note that $V_\beta(B)$ is
contained in the tubular neighborhood $V$ of $B$, for if
$d(h(x'),h(y'))<\beta$ with $h(x')\in B=h(B')$, then by the
quasi-isometry inequality $d'(x',y')<C_0\beta+D=\beta'$ so that
$y'\in V_{\beta'}(B')\subset V'$ and $h(y')\in V=h(V')$. Also
${\rm Vol}_{\beta}(B)\leq k=n_1k'$ since the image under $h$ of a
$\beta'$-ball on $L'$ is contained in a $C_0(\beta'+D)$-ball which
is covered by $n_1$ $\beta$-balls, and $k'$ $\beta'$-balls cover
$B'$. Thus $B$ satisfies the four conditions of Definition
\ref{BHP}.

By the bounded homology property for $L$, $C$ has Morse
$\beta$-volume $M(C,\beta)$ less than or equal to $K$, so that
there exists a Morse function $f: C\rightarrow [0,\infty)$ on $C$
whose level sets have $\beta$-volume at most $K$. Now $f\circ h$
is a Morse function on $C'$ whose level sets are taken onto those
of $f$ by $h$. Since $K$ balls of radius $\beta$ on $L$ suffice to
cover each level set of $f$, their images under $h^{-1}$ are
contained in $K$ balls of radius $\beta'$ on $L'$, and these balls
cover the corresponding level set of $f\circ h$. In other words,
$M(C',\beta')\leq K=K'$, as claimed. \qed

\section{Novikov's Theorem for embedded
vanishing cycles}\label{Novikov-gen} 

In this Section we prove the extension of Novikov's Theorem
(Theorem \ref{Nov-thm}, on the existence of Reeb components) and the
Corollary \ref{MorseReeb} which asserts that for every
$1$-parameter family of connected closed $(p-1)$-submanifolds
embedded in leaves in the interior of a Reeb component, each
bounding a compact region on its leaf, and for every
$\beta>0$, there is a common upper bound for the Morse
$\beta$-volumes of the regions that they bound, relative to any
fixed Riemannian metric on the Reeb component.

Throughout this section, we assume that
$\F$ is a $p$-dimensional topological foliation of a compact
$(p+1)$-dimensional manifold $M$ (so no differentiability is assumed), $B$ is a compact
connected $(p-1)$-dimensional manifold, and there is a foliated map
$$h: (B\times [0,1],\Hor)\rightarrow (M,\F)$$
where the horizontal and vertical foliations $\Hor$ and $\V$ of
$B\times [0,1]$ are given by the leaves $B\times \{t\}$ for $t\in
[0,1]$ and $\{x\}\times [0,1]$ for $x\in B$, respectively.
Furthermore we assume that $h_0: B\rightarrow L_0$ is an
embedding, where for all $t\in [0,1]$, $h_t: B\rightarrow L_t$ is
the map defined $h_t(b)=h(b,t)$ and $L_t$ is the leaf containing
$B_t=h_t(B)$.

The principal step in the Proof of Theorem \ref{Nov-thm}
is the following result. Let $C$ be a compact connected
$p$-dimensional manifold with boundary $B$.

\bpr If $h: (C\times (0,1]\cup B\times [0,1],\Hor,\V)\rightarrow (M,\F,\T)$
is a bifoliated embedding that cannot be extended over $C\times [0,1]$,
then the leaf $L_0$ containing $h(B\times \{0\})$ is the boundary
of a Reeb component whose interior is the union of the leaves meeting
$h(B\times (0,1])$.\label{Novlemma}\epr

\medskip

{\noindent\bf Proof of Theorem \ref{Nov-thm} from Proposition \ref{Novlemma}.}
In addition to the general hypotheses mentioned at the beginning of this section,
we suppose that one of the following conditions holds for every
$t>0$ sufficiently close to $0$ but not for $t=0$:
\begin{enumerate}
\item $B_t=h_t(B)$ is
the boundary of a compact $1$-connected region $C_t\subset L_t$;
\item $B_t=h_t(B)$ is
the boundary of a compact region $C_t\subset L_t$;
\item $\F$ and $B$ are oriented and $0=h_{t*}([B])\in H_{p-1}(L_t)$ (where $[B]$ is the
fundamental homology class of $B$); or
\item $0=h_{t*}([B])\in H_{p-1}(L_t;{\mathbb Z}_2)$ (where $[B]$ is the
fundamental homology class of $B$ with coefficients modulo $2$),
\end{enumerate}
and then we must show that the leaf $L_0$ is
the boundary of a Reeb component $R$ whose interior ${\rm Int}(R)$
is the union of the leaves $L_t$ for which $t>0$.

Let $\T$ be a $1$-dimensional foliation topologically transverse to $\F$.
We observe that it is possible to modify $h$ by a foliated homotopy,
moving $h(B\times \{t\})$ in the leaf $L_t$, with the homotopy fixed
on $h(B\times \{0\})$, so that the resulting mapping
$h: B\times [0,1],\Hor)\rightarrow (M,\F)$ restricts to a bifoliated mapping
$h_\epsilon: B\times [0,\epsilon],\Hor,\V)\rightarrow (M,\F,\T)$
for some $\epsilon>0$. Then it is easy to see that for some possibly smaller positive
$\epsilon$, $h_\epsilon$ will be an embedding. By reparametrizing the interval,
we may suppose that we have a bifoliated embedding
$h: B\times [0,1],\Hor,\V)\rightarrow (M,\F,\T)$. It is clear that any of the conditions
(1) through (4) that was satisfied by the original mapping $h$ will still hold
for the new bifoliated embedding $h$.

Next we observe that each of conditions (1), (3), and (4) imply condition (2).
Suppose (1), so that $B_t$ must be the boundary of a compact
$1$-connected region $C_t\subset L_t$ for every $t>0$. If $B_0$ bounded a
compact region $C_0$ on $L_0$, then by the Reeb Stability Theorem
$C_0$ would have a product foliated neighborhood and therefore would
be homeomorphic to $C_t$ with $t>0$ and hence $1$-connected, contrary to (1).
Thus (1) implies (2).

Condition (3) implies that for $t>0$, $B_t$
must be the boundary of a compact region $C_t$ contained in the leaf $L_t$,
since the $(p-1)$-dimensional cycle carried by $B_t$ is
a boundary on the $p$-dimensional leaf $L_t$. We are supposing in this case that
the foliation $\F$ is oriented, so if $B_0$ bounds a compact region $C_0$ in $L_0$,
then $C_0$ carries a relative homology class $[C_0]\in H_p(C_t,B_t)$ with
$\partial [C_0]=[B_0]$, thus making $[B_0]=0\in H_{p-1}(L_t)$, contrary to hypothesis.
This shows that (3) implies (2). The proof that (4) implies (2) is similar.
Thus we assume only condition (2), which includes the other three cases.

Now if for some values of $t$, $B_t$ bounds compact regions
on both sides in $L_t$, then $L_t$ will be a compact leaf. If there were
such values $t_n$ converging to $0$, then $L_0$ would also be a compact leaf
and for a sufficiently small $t_n$, the leaves $L_0$ and $L_{t_n}$ would bound
an $I$-bundle fibered by segments in the leaves of $\T$. Then $C_{t_n}$
would project along $\T$ onto a homeomorphic region $C_0\subset L_0$ so that
$\partial C_0=B_0$, contrary to hypothesis. Hence there is an
$\epsilon>0$ such that for each $t\in (0,\epsilon)$ $B_t$ bounds a compact region
on exactly one side in $L_t$.
Let $S_0$ be the set of $t\in (0,\epsilon)$ for which
$B_t$ bounds a compact region $C_t$ on the positive side in $L_t$,
and $S_1$ the set for which
$B_t$ bounds a compact region $C_t$ on the negative side, according to
a coherent transverse orientation of $B_t$ in $L_t$ which we choose arbitrarily. Then
$S_0\cap S_1=\emptyset$ and $S_0\cup S_1 = (0,\epsilon)$.

If there exists $0<\epsilon'\leq \epsilon$ such that $(0,\epsilon']
\subset S_0$ (or $(0,\epsilon']\subset S_1$),
then each of the corresponding compact regions $C_t$
must have trivial holonomy and consequently each $C_t$ will have a product neighborhood
foliated as a product. Joining these regions $C_t$ together we
obtain a mapping
$$h: (C\times (0,\epsilon']\cup B\times [0,\epsilon'],\Hor,\V)\rightarrow (M,\F,\T)$$
which does not extend over $C\times [0,\epsilon']$. Reparametrizing the interval
and applying Proposition \ref{Novlemma}
shows that $L_0$ is the boundary of a Reeb component with the leaves
$L_t$ for $t\in (0,\epsilon']$ in its interior, as claimed.
The remaining leaves containing $h(B\times \{t\}$ for the original mapping $h$ must also
be contained in the interior of the Reeb component, since the transverse orientation points
inwards along the boundary $L_0$.

In the remaining case, $0$ is a limit point of both $S_0$ and $S_1$, so
we can find a strictly decreasing sequence $t_n\searrow 0$ with
$t_n\in \overline S_0 \cap \overline S_1$. Suppose some $t_n\in S_0$
is a limit point of $S_1$ on the right. Considering the holonomy of
$C_{t_n}$, there must be an open interval $(t'_n,t'_n+\delta) \subset
S_1$ for some $\delta>0$ and some $t'_n\in S_0$ near to (and possibly
equal to) $t_n$. Applying the Proposition as before we find that
$L_{t'_n}$ bounds a Reeb component. The same conclusion holds if
$t_n$ is a limit point of $S_1$ on the left, and similarly if
$t_n\in S_1$. Thus we find a sequence $t'_n\searrow 0$ such that
every $L_{t'_n}$ bounds a Reeb component, which is impossible
since the boundary leaf of a Reeb component is compact and the
interior leaves are not. Hence this case does not occur, and the
Theorem is proven. \qed

\medskip

\noindent{\bf Proof of Proposition \ref{Novlemma}.} Suppose
that a mapping
$$h: (C\times (0,1]\cup B\times [0,1],\Hor,\V)\rightarrow (M,\F,\T)$$
is given, as in the Proposition. There must be a point
$x_0\in C\sm B$ such that $h$ does not extend to
the point $(x_0,0)$. (If
not, $h$ would extend uniquely and continuously over $C\times [0,1]$,
contrary to hypothesis.) Then the image of
the curve $t\mapsto h(x_0,t)$ lies on a leaf of $\T$, and as
$t\rightarrow 0$, the curve must have a point of accumulation $y_0$ on some
leaf $L$ in $M$. Hence we may find a strictly decreasing sequence
of numbers $t_n\searrow 0$ such that the sequence
$y_n=h(x_0,t_n)$ converges to $y_0$ as $n\rightarrow\infty$. Let
$V$ be a connected open neighborhood of $y_0$ on its leaf $L$.

\blm It is possible to choose the sequence $\{t_n\}$, the limit
point $y_0$, its neighborhood $V\subset L$, and $\epsilon>0$, so that
\begin{enumerate}
\item the leaf $L$ containing $y_0$ is distinct from $L_0$;
\item for every $n$, $t_n<\epsilon$ and $y_n=h(x_0,t_n)\in V$; and
\item $V$ is disjoint from $h(B\times [0,\epsilon])$.
\end{enumerate}
\elm

\noindent{\bf Proof.} If it happens that $y_0\in L_0$, by a small
change of the values of the numbers $t_n$, we may move the points
$y_n=h(x_0,t_n)$ along $\T$ so that the sequence $y_n$ converges to another point (still
denoted $y_0$) on another leaf $L$. Then since $y_0\notin h(B\times \{0\})$, for a sufficiently
small $\epsilon>0$, $y_0$ will not lie on the set $h(B\times [0,\epsilon])$,
and we may choose an open connected neighborhood $V$ of $y_0$ on
its leaf $L$ whose closure $\bar V$ is disjoint from the compact set $h(B\times
[0,\epsilon])$. Then, slightly changing the values of the numbers
$t_n$ and possibly passing to a subsequence, we may guarantee that
$t_n<\epsilon$ and $y_n=h(x_0,t_n)\in V$. \qed

\medskip

Now $y_n$ lies on the intersection of $V$ with the region
$C(n)=h(C\times \{t_n\})$ on $L$, and $V$ is disjoint from its
boundary $B(n)=h(B\times \{t_n\})=\partial C(n)$, so
$V\subset C(n)$. The connected submanifolds $B(m)=\partial C(m)$)
are pairwise disjoint and $y_0\in C(n)\cap C(m)$, so either $C(n)\subset
C(m)$ or $C(m)\subset C(n)$.  The sets $B(m)$ are separated by a
positive distance on the leaf $L$, so only finitely many of the
$C(m)$ can be contained in the compact set $C(n)$. Thus for each
$n$ there exists some $n'>n$ such that $C(n)\subset C(n')$.

\blm The leaf $L_0$ is compact. \elm

\noindent{\bf Proof.} If $L_0$ is not compact, then there is
a simple closed curve $\gamma$ transverse to $\F$ in the positive
direction that intersects $L_0$ in a single point $\gamma(0)=z$;
we can choose $z$ to be near to but not on the set $B_0$, and on
the side of $B_0$ in $L_0$ on which each submanifold $B_t$ bounds
$C_t$. We may isotope $\gamma$ slightly so that for some small
positive numbers $\epsilon_1$ and $\epsilon_2\leq\epsilon$, the
segment $\gamma([0,\epsilon_1])$ lies on a leaf of $\T$ and
$\gamma$ is disjoint from $h(B\times [0,\epsilon_2])$. Then for a
sufficiently large index $n$, with $t_{n'}<t_n<\epsilon_2$,
$\gamma$ will enter into the immersed region $h(C\times
[t_{n'},t_n])$ at a point $\gamma(s)\in C(n')$ for some $s\in
(0,\epsilon_1)$. Now $\gamma$ can never exit from that region,
whose boundary is contained in $C(n')\sm {\rm Int}\ C(n)\cup
h(B\times [t_{n'},t_n])$, for $\gamma$ cannot exit along $C(n')\sm
C(n)$ where the transverse orientation enters, and $\gamma$ is
disjoint from $h(B\times [t_{n'},t_n])\subset
h(B,[0,\epsilon_2])$. This contradiction shows that $L_0$ must be
compact. \qed

\medskip

Let $N$ be a positive one-sided tubular neighborhood of $L_0$
fibered by segments in leaves of $\T$ with projection map
$p:N\rightarrow L_0$. Let $L_0'\subset N'$ be the result of
cutting $L_0$ along $B_0$ and cutting $N$ along $p^{-1}(B_0)$ to
get $N'$, so that $\partial L_0'=B_0^+\cup B_0^-$, two disjoint
copies of $B_0$, and $N'$ is a positive one-sided tubular
neighborhood of $L_0'$. The positive holonomy of the compact leaf
$L_0'$ must be trivial, for the leaves near to $L_0'$ are
contained in the compact sets $C_t$ and thus are compact. Hence a
smaller compact tubular neighborhood $N_0'\subset N'$ will be
foliated as a product, say by leaves $D_t$ with boundary $\partial
D_t=B_t^+\cup B_{f(t)}^-$, where $f$ is a function defined on a
small positive one-sided neighborhood of $0$ in $[0,\epsilon_2)$
and $p(B_s^\pm)=B_0^\pm$ for every sufficiently small $s$. For
definiteness we choose the notation so that $f(t)<t$ and
consequently $C_t\subset C_{f(t)}$. This holds for all $t$ less
than or equal to some $t_0\leq \epsilon_2$. As $t$ varies in
$[0,t_0]$ there is defined a continuous function
$f:[0,t_0]\rightarrow [0,t_0]$ such that $f(t)<t$ for every $t>0$,
while in the limit $f(0)=0$. Then it is clear that $C_t\cup D_t
=C_{f(t)}$ with $C_t\cap D_t= B_t$.

Observe that the inclusion $i_t: C_t\subset C_{f(t)}$ defines an
embedding $\phi= h_{f(t)}^{-1}\circ i_t\circ h_t: C\rightarrow C$ that
does not depend on $t\in (0,t_0]$, since moving along the leaves
of $\T$ produces the same result for each $t$. Thus $h: C\times
(0,t_0]\rightarrow M$ passes to a quotient immersion $\bar h: R_0
\rightarrow M$, where $R_0= C\times (0,t_0]/\{(x,t)\sim
(\phi(x),f(t))\}$. For each $t\in (0,t_0]$, let $L(t)$ be the
image of the set $\cup_{n=0}^\infty\ C\times \{f^n(t)\}$ in $R_0$.
Then $\bar h|_{L(t)}: L(t)\rightarrow L_t$ is a homeomorphism, for the
identifications correspond to the inclusions $C_t\subset
C_{f(t)}$; there cannot be any further identifications, since the
regions $C_t$ were chosen to be embedded in leaves of $\F$ in $M$,
and the union of the sets $C_{f^n(t)}$ exhausts the leaf $L_t$. No
region $C_s$ with $s>t$ that is not one of the sets $C_{f^n(t)}$
can meet the leaf $L_t$, for then
there would be an index $n\geq 0$ such that $C_{f^n(t)}\subset C_s
\subset C_{f^{n+1}(t)}$, which is impossible since $D_{f^n(t)}$
contains no $B_s$ in its interior. Consequently $\bar h:
R_0\rightarrow \cup\ \{L_t|\ t\in (0,t_0)\}$ is a bijection, which is
easily seen to be a homeomorphism.

The map $p: R_0 \rightarrow S^1=[t,f(t)]/\{t\sim f(t)\}$, defined
by setting $p(L_s)=s$ if $s\in [t,f(t)]$, is well defined, and it is a
fibration whose local product structure is given by translations
by holonomy mappings along leaves of $\T$. A small positive
compact tubular neighborhood $N_0$ of $L_0$ will meet $R_0$ in
$N_0\sm L_0$. Hence $R=N_0\cup R_0$ is compact, since it coincides
with the union of the two compact sets $N_0$ and $h(C\times
[t_1,f(t_1)])$ for some sufficiently small positive $t_1$;
furthermore $R=L_0\cup R_0$ and $\partial R=L_0$. Thus we have
shown that $R$ is a compact manifold with boundary $L_0$, and its
interior $R_0$ fibers over the circle. Since the boundary is
connected, there exists a transverse orientation pointing inwards,
so $R$ is a Reeb component. Finally, all the leaves $L_t$ for
$t\in (0,1]$ in the original parametrization of the interval
are contained in $R_0$ because the transverse
orientation points inwards along $\partial R$; for any point $z\in
B$ the curve $t\mapsto h(z,t)$ lies in $R_0$ for small values of
$t$ and as $t$ increases it must be entirely contained in $R_0$.
\qed

\medskip
\begin{figure}[H]\label{phionC}
\centering
\includegraphics*[width=.4\linewidth]{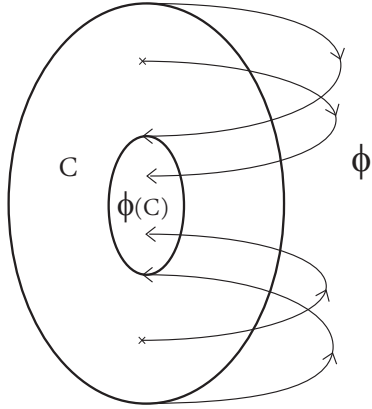}
\caption{$\phi:C\hookrightarrow \Int\ C$.}
\end{figure}
Note that in the preceding proof the Reeb component $R$ shown to
exist was obtained from the map $\phi: C\rightarrow {\rm Int}\ C$
and the contraction $f: [0,t_0] \rightarrow [0,t_0]$, by a
construction which we shall now describe. It is not difficult to
show that such a construction gives {\bf all} Reeb components, up
to foliated homeomorphism (see \cite{AHS} for the proof), although
we shall not use that fact here.

\medskip \noindent {\bf Construction of Reeb components.}
Let $C$ be a connected compact $p$-manifold $C$ with connected
boundary $\partial C=B$ and let $\phi: C\rightarrow {\rm Int}\ C$
be an embedding of $C$ into its own interior. (See Figure 8.) This
construction generalizes to the case in which $B$ is not
connected, but an extra condition is required, and we do not need
this generalization here, so we omit it. Consider the product
$C\times [0,1]$ with the product foliation whose leaves are
$C\times \{t\}$.

\noindent On the submanifold
$$C'= C\times [0,1]\ \sm\ \phi({\rm Int}\ C)\times\{0\}$$ let $\sim$
be the equivalence relation generated by setting $(x,s)\sim
(\phi(x),f(s))$ for every $(\phi(x),f(s))\in C'$, where $f:
[0,1]\rightarrow [0,1]$ is a continuous embedding such that
$f(0)=0$ and $f(s)<s$ for every $s>0$. The compact set
$$(C\sm
\phi({\rm Int}\ C))\times [0,1]\ \cup\ C\times [f(1),1]$$

\begin{figure}[H]\label{Cwithidents}
\centering
\includegraphics*[width=.4\linewidth]{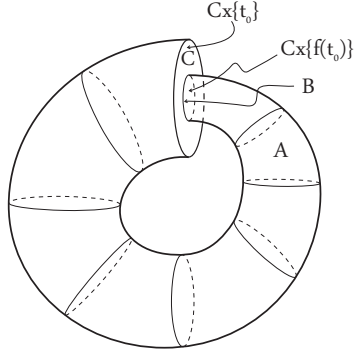}
\caption{$C\times [f(t_0), t_0]$ with identifications.}
\end{figure}
\begin{figure}[H]\label{otherview}
\centering
\includegraphics*[width=.3\linewidth]{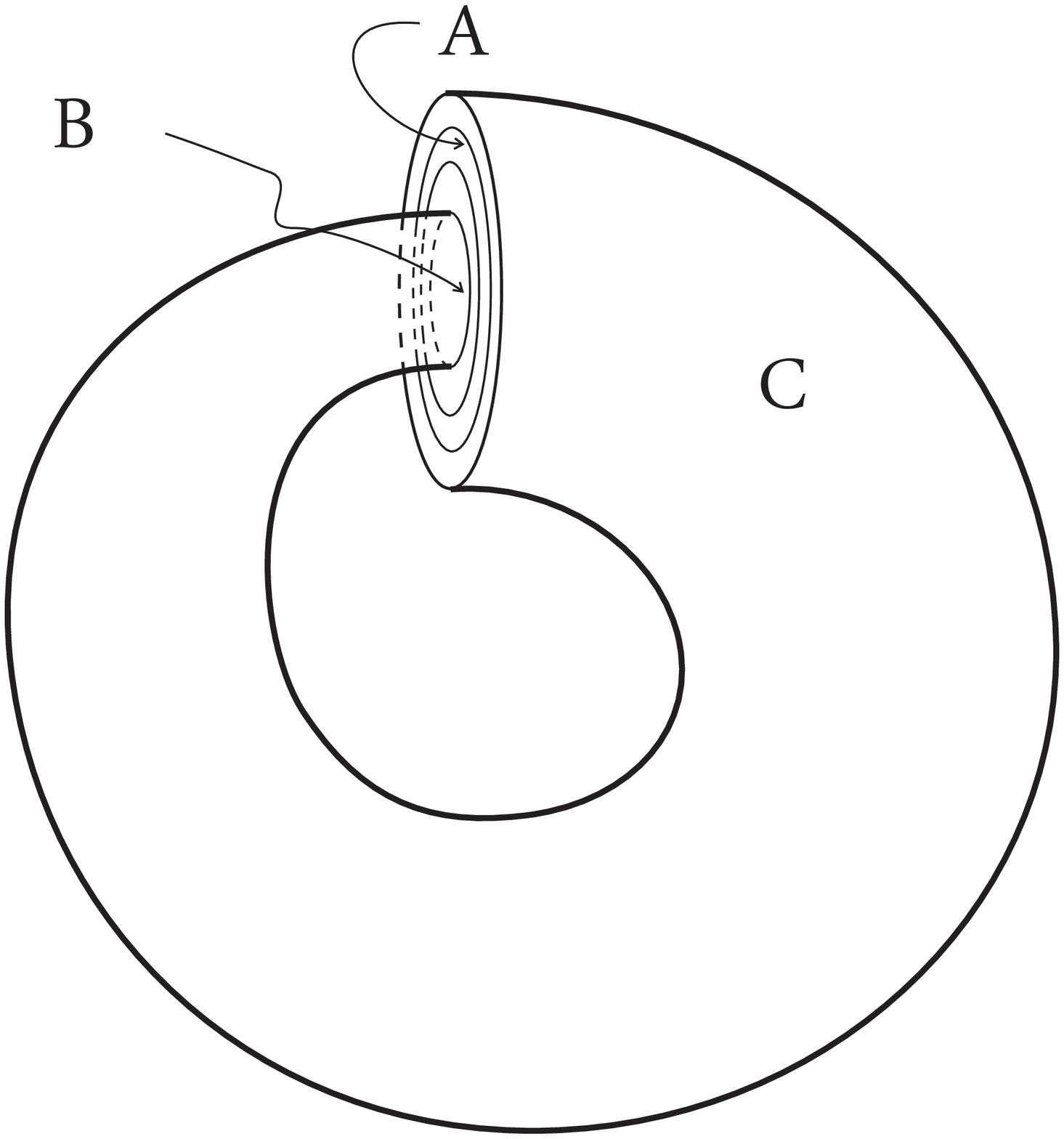}
\caption{A diffeomorphic image of the same set showing part of the
Reeb component.}
\end{figure}
\noindent projects onto the quotient $R=C'/\sim$, so $R$ is
compact. It is not difficult to check that $R$ is a compact
$(p+1)$-manifold with boundary endowed with a codimension one
foliation $\R$ induced by the horizontal foliation on $C'$, and
that $(R,\R)$ is a Reeb component. Figure 9 shows part of the Reeb
component, the image of $C\times [f(t_0),t_0]$ with $C\times
\{t_0\}$ identified with $\phi(C)\times \{f(t_0)\}$ by the
equivalence relation $\sim$, foliated by the images of the sets
$C\times \{t\}$. Figure 10 shows a diffeomorphic image of the same
compact region with corners, but this view suggests how the Reeb
component is built up as $t_0$ decreases to $0$. The set $(C \sm
\phi(\Int\ C))\times \{0\}$ projects onto the boundary of the Reeb
component. (We remark that the foliated homeomorphism type is
independent of the choice of the embedding $f$ since $f$ is
topologically conjugate to any other embedding with the same
properties.)

Corollary \ref{MorseReeb} will follow from the following result,
since the Reeb component constructed in the Proof of Theorem
\ref{Nov-thm} satisfies its hypotheses. Note that we are assuming
that $\F$ is $C^{2,0}$.

\bpr Let $C$ be a compact connected $p$-manifold with connected
boundary and let $\phi: C\rightarrow \Int(C)$ be a $C^2$ embedding
into its interior. Then for every $\beta>0$ there is a constant
$K>0$ such that $M(C_t,\beta)\leq K$ for every region $C_t$
appearing in the above construction of the Reeb component
$(R,\R)$. \epr

Corollary \ref{MorseReeb} will follow from the following result, since the
Reeb component constructed in the Proof of Theorem \ref{Nov-thm} satisfies
its hypotheses. Note that we are assuming that $\F$ is $C^{2,0}$.

\bpr Let $C$ be a compact connected $p$-manifold with connected boundary and
let $\phi: C\rightarrow \Int(C)$ be a $C^2$ embedding into its interior.
Then for every $\beta>0$ there is a constant $K>0$ such that $M(C_t,\beta)\leq K$
for every region $C_t$ appearing in the
above construction of the Reeb component $(R,\R)$. \epr

\medskip

\noindent {\bf Proof}. Let $(R,\R)$  be the Reeb component
constructed from $\phi: C\rightarrow \Int(C)$, as above, where $C$
and $\phi$ are smooth of class $C^2$. The boundary $L_0=\partial
R$ is diffeomorphic to $D/\{x \sim \phi(x)\}$, where $D= C\sm
\phi({\rm Int} \ C)$. Let $B_0$ be the image of $\partial
D=\partial C \cup\partial \phi(C)$ in $L_0$ under the
identification. Clearly $B_0$ is two-sided in $L_0$, so by using a
tubular neighborhood of $B_0$ in $L_0$ we may find a smooth map
$g_0: L_0\rightarrow S^1$ with $1\in S^1$ as a regular value and
such that $g_0^{-1}(1)=B_0$. By a small perturbation we may
suppose that $g_0$ is a smooth Morse function. (For convenience,
we first consider Morse functions with values in $S^1$ rather than
in ${\mathbb R}$.) Since $L_0$ is compact, there is an upper bound
$K_1$ on the $\beta$-volumes ${\rm Vol}_{\beta}(g_0^{-1}(z))$ for
all $z\in S^1$. Extend $g_0$ to $g:N\rightarrow S^1$, where $N$ is
a small tubular neighborhood of $L_0$ in $R$, by setting
$g=g_0\circ p$ where $p:N\rightarrow L_0$ is the projection along
leaves of the transverse foliation $\T$ induced by the vertical
foliation on $C'= C\times [0,1]\ \sm\ \phi({\rm Int}\
C)\times\{0\}$. Let $\tilde N$ be the cyclic cover of $N$
corresponding to the map $g_0:N\rightarrow S^1={\mathbb
R}/{\mathbb Z}$ and let $\tilde g: \tilde N\rightarrow {\mathbb
R}$ be the natural lift of $g$ to the cyclic cover. We let $C_t$,
$B_t$, and $D_t$ be the images in $R$ of $C\times \{t\}$,
$\partial C\times \{t\}$, and $(C\sm \partial C)\times \{f(t)\}$
under the identification. As in the proof of Theorem
\ref{Novikov-gen}, $C_t \cup D_t = C_{f(t)}$ with $C_t \cap D_t =
B_t$ and $\partial D_t = B_t \cup B_{f(t)}$. For each $t$, the
region $C_t\cap N$ can be lifted to $\tilde N$ by lifting $D_t,
D_{f(t)}$, etc., successively. For a sufficiently small tubular
neighborhood $N$ of $L_0$, $\tilde g$ will restrict to a Morse
function $\tilde g_t$ on $C_t\cap N$, and $K_1$ will be an upper
bound for the $\beta$-volume ${\rm Vol}_{\beta}(\tilde
g_t^{-1}(r))$ of each level set $\tilde g_t^{-1}(r)$ for $r\in
{\mathbb R}$, since $\tilde g_t^{-1}(r)$ is a compact set close to
$g_0^{-1}(r\ {\rm mod}\ 1)$ which is covered by at most $K_1$ open
balls of radius $\beta$. Extend $\tilde g_t$ to a Morse function
$\hat g_t:C_t\rightarrow {\mathbb R}$.

Next, let $S_s$ be the image of $C \times [f(s),s]$ in ${\rm
Int}(R)$. The sets $S_s$ are nested and their interiors cover the
compact set $R\sm {\rm Int}\ N$, so for a sufficiently small $s>0$
the union ${\rm Int}\ N\cup {\rm Int}\ S_s$ will be the whole
manifold $\bar R$. The images $C_t$ of the sets $C \times \{t\}$
in the leaf $L_t$ of $\R$ are compact and vary continuously, so
there is a common upper bound $K_2$ for their $\beta$-volumes for
all $t\in [f(s),s]$. Finally each level set $\hat g_t^{-1}(r)$ for
$r\in {\mathbb R}$ is contained in the union of $\tilde
g_t^{-1}(r)\subset \tilde N$ and $S_t$, so its $\beta$-volume is
at most $K_1+K_2$, a common upper bound for the $\beta$-volumes of
the level sets $\hat g_t^{-1}(r)$ for all $t$ and $r$. It
follows that $K_1+K_2$ is a common upper bound for the Morse
volumes of the sets $C_t$, as claimed. \qed

\end{document}